\documentclass[12pt,reqno]{amsart}

\usepackage[T1]{fontenc}
\usepackage[english]{babel}
 \usepackage{colonequals} 

\usepackage[margin=1.3in,footskip=.2in]{geometry}
\usepackage{amsmath,amssymb,amscd,amsfonts,mathtools,mathrsfs,tikz,caption,float,tikz-cd}
\usetikzlibrary{arrows,patterns,calc}
\usetikzlibrary{decorations.markings}
\usepackage[pdffitwindow=true,pagebackref=true,linkcolor=blue,colorlinks,citecolor=blue,urlcolor=blue]{hyperref}

\tikzstyle{arrowmid}[0.5]=[decoration=
{markings,mark=at position #1 with {\arrow{>}}},
postaction={decorate}]

\newtheorem{theorem}{Theorem}[section]
\newtheorem{lemma}[theorem]{Lemma}
\newtheorem{definition}[theorem]{Definition}
\newtheorem{proposition}[theorem]{Proposition}
\newtheorem{corollary}[theorem]{Corollary}

\newcommand{\figref}[1]{\hyperref[#1]{Figure \ref{#1}}}
\newcommand{\lemref}[1]{\hyperref[#1]{Lemma \ref{#1}}}
\newcommand{\thmref}[1]{\hyperref[#1]{Theorem \ref{#1}}}
\newcommand{\conjref}[1]{\hyperref[#1]{Conjecture \ref{#1}}}
\newcommand{\propref}[1]{\hyperref[#1]{Proposition \ref{#1}}}
\newcommand{\cororef}[1]{\hyperref[#1]{Corollary \ref{#1}}}
\newcommand{\defref}[1]{\hyperref[#1]{Definition \ref{#1}}}
\newcommand{\rmkref}[1]{\hyperref[#1]{Remark \ref{#1}}}
\newcommand{\qref}[1]{\hyperref[#1]{Question \ref{#1}}}
\newcommand{\secref}[1]{\hyperref[#1]{\S\ref{#1}}}
\newcommand{\appref}[1]{\hyperref[#1]{Appendix \ref{#1}}}

\newcommand{\scal}[1]{\left\langle#1\right\rangle}

\newcommand{\op}[1]{\left(#1\right)}

\newcommand{\R}{\mathbb{R}}
\newcommand{\C}{\mathbb{C}}

\newcommand{\h}{\mathbb{H}}
\newcommand{\Q}{\mathbb{Q}}
\newcommand{\Z}{\mathbb{Z}}
\newcommand{\N}{\mathbb{N}}

\newcommand{\eps}{\varepsilon}
\newcommand{\im}{\mathrm{Im}}
\newcommand{\re}{\mathrm{Re}}

\newcommand{\Li}{\mathrm{Li}}

\newcommand{\bbm}{\begin{bmatrix}}
\newcommand{\ebm}{\end{bmatrix}}
\newcommand{\bpm}{\begin{pmatrix}}
\newcommand{\epm}{\end{pmatrix}}
\newcommand{\bsm}{\left(\begin{smallmatrix}}
\newcommand{\esm}{\end{smallmatrix}\right)}
\newcommand{\bsbm}{\left[\begin{smallmatrix}}
\newcommand{\esbm}{\end{smallmatrix}\right]}

\newcommand{\fa}{\frak{a}}
\newcommand{\fb}{\frak{b}}

\newcommand{\gs}{\sigma}

\newcommand{\G}{\Gamma}
\newcommand{\g}{\gamma}
\newcommand{\bk}{\backslash}

\newcommand{\cF}{\mathcal{F}}
\newcommand{\cH}{\mathcal{H}}
\newcommand{\cS}{\mathcal{S}}

\newcommand{\sign}{\mathrm{sign}}

\newcommand{\tr}{\mathrm{tr}}
\newcommand{\1}{\mathrm{1}}
\newcommand{\SL}{\operatorname{SL}}
\newcommand{\GL}{\operatorname{GL}}

\newcommand{\PSL}{\operatorname{PSL}}

\begin{document}

% Enter full title and short title for running headers
\title{Windings of Prime Geodesics}

\author{Claire Burrin}
\address{Institute of Mathematics, University of Zurich, Switzerland}
\email{claire.burrin@math.uzh.ch}

\author{Flemming von Essen}
\email{flemmingvonessen@gmail.com}

%\date{\today}

\begin{abstract}
The winding of a closed oriented geodesic around the cusp of the modular orbifold is computed by the Rademacher symbol, a classical function from the theory of modular forms. In this article, we introduce a new construction of winding numbers to record the winding of closed oriented geodesics about a prescribed cusp of a general cusped hyperbolic orbifold. For various arithmetic families of surfaces, this winding number can again be expressed by a Rademacher symbol, and access to the spectral theory of automorphic forms yields statistical results on the distribution of  closed (primitive) oriented geodesics with respect to their winding.
\end{abstract}

\maketitle

%\item[(iii)] Note the use of \verb"\affil" and \verb"\affilnum" to link names and
%addresses. The address and/or email address of the author for
%correspondence is defined by \verb"\correspdetails".

\section{Introduction}

The winding of a closed oriented geodesic around the cusp of the modular orbifold $\PSL_2(\Z)\bk\h$ can be expressed by an explicit function $\Psi$ from the theory of modular forms called the Rademacher symbol. In Figure \ref{fig} below, the coloured path represents a portion of the dashed hyperbolic geodesic folded onto the standard fundamental domain for the modular group $\PSL_2(\Z)$.

\begin{figure}[ht]
\centering
\includegraphics[scale=0.6]{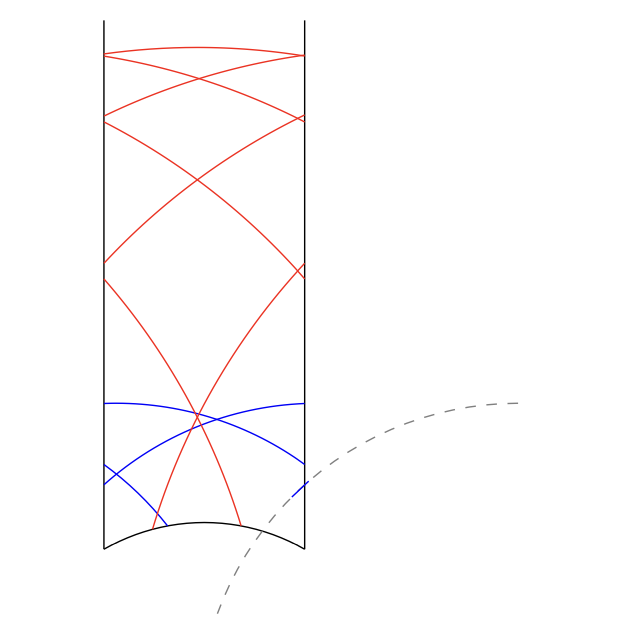}
%\begin{tikzpicture}[scale=2]
%% real axis
%%\draw (-1,0) -- (1,0);
%% fundamental domain
%\draw [domain=pi/3:2*pi/3] plot ({cos(\x r)},{sin(\x r)});
%\draw [domain={sin(pi/3 r)}:3.5] plot ({cos(pi/3 r)},\x);
%\draw [domain={sin(2*pi/3 r)}:3.5] plot ({cos(2*pi/3 r)},\x);
%
%% hyperbolic geodesic
%\draw[gray,dashed] [domain=90:160] plot ({50/32+(51/32)*cos(\x)},{(51/32)*sin(\x)});
%
%% first iteration of winding
%\draw[blue] (-0.5,1.1879) arc ({atan2(1.1879,-17/16)} : {atan2(1.5925,-1/16)} : 51/32) ;
%\draw[blue] (-0.5,1.5925) arc ({atan2(1.5925,-1/16)} : {atan2(1.2888,15/16)} : 51/32) ;
%\draw[blue] (-0.5,1.2888) arc ({atan2(1.2888,15/16)} : {atan2(0.9831,1.2544)} : 51/32) ;
%
%% second iteration of winding
%\draw[red] (0.18308,0.9831) arc ({atan2(0.9831,3.218)} : {atan2(2.21278,2.535)} : 1632/485) ;
%\draw[red] (0.5,2.21278) arc ({atan2(2.21278,2.53505)} : {atan2(2.994,1.53505)} : 1632/485) ;
%\draw[red] (0.5,2.994) arc ({atan2(2.994,1.53505)} : {atan2(3.32216,0.535050)} : 1632/485) ;
%\draw[red] (0.5,3.32216) arc ({atan2(3.32216,0.535050)} : {atan2(3.32216,-0.465)} : 1632/485) ;
%\draw[red] (0.5,3.3269) arc ({atan2(3.3269,-0.465)} : {atan2(3.029,-1.465)} : 1632/485) ;
%\draw[red] (0.5,3.029) arc ({atan2(3.029,-1.465)} : {atan2(2.2906,-2.465)} : 1632/485) ;
%\draw[red] (0.5,2.2906) arc ({atan2(2.2906,-2.465)} : {atan2(0.966,-3.2235)} : 1632/485) ;
%
%% base point 
%\draw[blue] ({50/32+(51/32)*cos(3*pi/4 r)},{(51/32)*sin(3*pi/4 r)}) -- (0.52,1.2058) ;
%%\draw[blue] ({50/32+(51/32)*cos(3*pi/4 r)},{(51/32)*sin(3*pi/4 r)}) node[left] {$x$} ;
%\end{tikzpicture}
\captionof{figure}{}\label{fig}
\end{figure}

Travelling left to right along the blue segment, we record $a_1=3$ windings around the cusp until eventually hitting the lower boundary arc of the fundamental domain. The geodesic trajectory then continues along the red segment,  where we now go from right to left, and we  record the $a_2=7$ next windings until again reaching the lower boundary arc, and so on. If the dashed geodesic is the axis of a closed geodesic $C$, this process yields a finite sequence $(a_1,\dots,a_n)$ that encodes the {\em oriented} winding around the cusp. For the minimal even length such sequence, the Rademacher symbol $\Psi(C)$ is equal to the alternating sum 
$$
\Psi(C) = a_1-a_2 +a_3 -a_4 +\dots+ a_n.
$$
This observation --- the (oriented) winding of a closed geodesic about the cusp is realized by the Rademacher symbol --- allows to deduce counting and distribution results for the windings of modular geodesics. This was implemented by Sarnak \cite{Sarnak2010} (see  \cite{Mozzochi2013} for details) who obtained the following result, among others that we will discuss later in this introduction. Let $\Pi(T)$ be the countable set of prime (i.e., primitive, closed, and oriented) geodesics on $\PSL_2(\Z)\bk \h$ of length up to $T$.  For $n\in\Z$ fixed, the asymptotic density of prime geodesics with winding number $\Psi(C) =n$ is 
$$
\frac{\pi_n(T)}{\pi(T)}\colonequals\frac{\#\{C\in\Pi(T): \Psi(C)=n\}}{\#\Pi(T)}\, \sim\, \frac{1}{3} \frac{T}{T^2 +\op{\tfrac{\pi n}{3}}^2}
$$
as $T\to\infty$. Together with the prime geodesic theorem, this shows that there are infinitely many prime geodesics winding $n$ times for each choice of $n\in\Z$. The distribution of the winding numbers is symmetric about $n=0$ --- reflecting that the sign of the winding number changes if we reverse the direction of the geodesic --- and peaks at $n=0$. We will soon see that this statistical behavior persists for a natural extension of this winding number to more general hyperbolic surfaces.

\vspace{.2cm}
We recall the following notion of winding for closed curves on a cusped hyperbolic surface $M$, following \cite{Reinhart1960,Chillingworth1972}. Fix a non-vanishing continuous vector field $X$ on $M$.\footnote{ The restriction to noncompact surfaces is essential for this construction; compact surfaces may not admit such vector fields. For an algebraic construction of winding numbers for compact surfaces, see \cite{Huber2012}.} Then for each closed parametrized curve $C(t)=(c(t),\theta(t))$ in the unit tangent bundle ${\rm T}^1 M$, given by position $c(t)$ and direction $\theta(t)$ at time $t$, the winding number of $C$ relative to the vector field $X$ is defined to be the total variation of the angle function 
$$
\theta(t)\colonequals \measuredangle(\theta(t),X(c(t)))
$$ 
over one revolution. This winding number is independent of the choice of parametrization and classifies nontrivial homotopy classes of closed curves in ${\rm T}^1 M$ \cite{Chillingworth1972,Smale1958}. %Adjust to reality
The main object of this paper is to introduce a similar winding number function with the vector field $X$ replaced by an automorphic form $f$ that would generalize the natural winding number function provided by the Rademacher symbol for the modular orbifold.

\vspace{.2cm}
Our construction goes as follows. Let $M=\G\bk \h$ be a cusped hyperbolic surface. 
For simplicity, we fix a fundamental domain with a cusp at $\infty$.\footnote{ We may assume the existence of such a cusp up to replacing $M$ by an isometric surface.} An automorphic form of (integer) weight $2k$ for the action of $\G$ is a function $f:\h\to\C$ that is holomorphic on $\h$ and at the cusps, and such that $f(\g z)(\g' z)^k =f(z)$ for all $\g\in\G$, $z\in\h$. The (finite-dimensional) space of such forms is denoted $M_{2k}(\G)$. In addition to these stringent conditions, we further require that $f\in M_{2k}(\G)$ be
\begin{itemize}
\item[(i)] nowhere-vanishing, i.e., $f(z)\neq0$ for all $z\in\h$;
\item[(ii)]  and $f(x+iy)\to0$ as $y\to\infty$. 
\end{itemize}
By a standard procedure we may lift $f$ to a (unique) continuous function $F:{\rm T}^1M\to \C^*$, where $\C^*=\C\setminus\{0\}$ denotes the punctured plane. Then each closed curve $C$ in ${\rm T}^1 M$ is mapped to a continuous oriented closed curve $F(C)\subset\C$. By construction, the curve $F(C)$ never passes through the origin of the plane, which corresponds to the cusp at $\infty$. We define the {\bf winding of $C$ relative to $f$} to be the topological index of the planar curve $F(C)$ about the origin, i.e.,
$$
{\rm ind}(F(C)) = \frac{1}{2\pi i}\int_{F(C)} \frac{dz}{z}.
$$
This integer is independent of parametrization and homotopy-invariant. A feature of negative curvature is that each homotopy class contains a unique geodesic representative and we henceforth restrict $C$ to the family $\Pi$ of prime geodesics on $M$. Further, given two forms $f,g\in M_{2k}(\G)$, the associated winding numbers ${\rm ind}(F(C))$ and ${\rm ind}(G(C))$ differ by $h([c])$, where $h$ is a homomorphism on the homology group ${\rm H}^1(M)$. This observation made, our goal is to construct an explicit automorphic form\footnote{ We will eventually see in the course of the constructive proof that holomorphic can be relaxed to real-analytic.} $f$ with properties (1)-(2) such that the arising winding number ${\rm ind}(F(C))$ coincides with a natural extension of the Rademacher symbol $\Psi$. 

There has recently been a renewed interest in variations around the Rademacher symbol, e.g., the linking number introduced by Duke, Imamoglu, and Toth in \cite{DukeImamogluToth2016}, 
 the Dedekind--Rademacher cocycle studied by Darmon, Pozzi, and Vonk in \cite{DarmonVonk2021}, or the extension of Rademacher symbols to Fuchsian groups introduced by the first named author \cite{Burrin2022} in relation to the Manin--Drinfeld theorem, and studied by Matsusaka and Ueki \cite{MatsusakaUeki2023} in the special case of triangle groups from the viewpoint of linking numbers. From \cite{Burrin2022} we know that each cusp $\fa$ of $M$ gives rise to a generalization $\Psi_\fa$ of the Rademacher symbol. 
 We only consider the Rademacher symbol $\Psi:=\Psi_\infty$ associated to the cusp at $\infty$ --- for the case of a general cusp, the reader is referred to \cite{Burrin2022}. For the purpose of this introduction, we only state a simplified version of our main result, and postpone the full statement to Section \ref{discu} (see Theorem \ref{Thm1b}) as well as a discussion on the existence and construction of the required forms.

\begin{theorem}\label{Thm0}
Let $M=\G\bk\h$ be a cusped hyperbolic orbifold with finite hyperbolic area $V={\rm area}(M)$ and let $\Psi$ be its associated Rademacher symbol, as above. There exists an automorphic form $f$, possibly nonholomorphic, of even weight $k\in2\N$ that is nowhere-vanishing, for which $f(z)\to0$ as $y\to\infty$, and such that the winding number relative to $f$ is either  
\begin{align}\label{bullet}
{\rm ind}(F(C)) = \frac{kV}{4\pi}\Psi(C)
\end{align}
or a linear combination of values of $\tfrac{kV}{4\pi}\Psi$.
\end{theorem}

For the modular orbifold $\PSL_2(\Z)\bk\h$ of area $\tfrac{\pi}{3}$, the modular form $f$ coincides with the modular discriminant $\triangle$ of weight $k=12$, and we recover the classical Rademacher symbol $\Psi$. Further, we will see that we are in the situation of \eqref{bullet} whenever $\G\in\mathcal{G}$, where 
 $\mathcal{G}$ is the collection of the following families of Fuchsian groups: congruence subgroups of the modular group, maximal arithmetic noncocompact Fuchsian groups (under inclusion), genus 0 noncocompact Fuchsian groups (e.g., Hecke triangle groups). In these situations we are able to extend Sarnak's strategy for counting and we obtain the following results. 

 \begin{theorem}\label{Thm1}
 Let $\G\in\mathcal{G}$ and consider the winding number given by \eqref{bullet}. Recall that $\pi_n(T)$ is the number of prime geodesics in $\Pi(T)$ with winding $n$. There exists $\delta\in(0,\tfrac12]$ such that 
$$
\pi_n(T) = \frac{4}{k T}\int_2^{e^T} \frac{\log t}{(\log t)^2 +\op{\tfrac{4\pi n}{k}}^2}\, dt\, + O\op{\frac{e^{T(1-\delta/2)}}{T}}
$$
as $T\to\infty$. The implied constant does not depend on $n$.
\end{theorem}

\begin{corollary}
 Let $\G\in\mathcal{G}$ and consider the winding number given by \eqref{bullet}. For $n\in\Z$ fixed, the asymptotic density of prime geodesics with $n$ windings is
 $$
 \frac{\pi_n(T)}{\pi(T)}\, \sim\,  \frac{4}{kT} \frac{T}{T^2+\op{\tfrac{4\pi n}{k}}^2}
 $$
 as $T\to\infty$.
 \end{corollary}

\begin{theorem}\label{Thm:Cauchy}
Let $\G\in\mathcal{G}$ and consider the winding number given by \eqref{bullet}. The limiting distribution of the ratio of winding-to-length for prime geodesics is Cauchy. More precisely, for any interval $[a,b]$ of the real line we have
$$
\lim_{T\to\infty} \frac{\#\{ C\in \Pi(T) : \tfrac{k}{4\pi}a \leq \tfrac{{\rm ind}(F(C))}{\ell_C} \leq \tfrac{k}{4\pi}b\}}{\pi(T)}
= \int_{a}^{b} \frac{du}{\pi(1+u^2)}.
$$
\end{theorem}

\vspace{.2cm}
These results generalize \cite[Theorems 1 and 3]{Sarnak2010}, which are stated for $\G=\PSL_2(\Z)$, $\delta=\tfrac12$, and $k=12$. We believe that the presence of cusps accounts for the fat tails of the Cauchy distribution; given a long closed geodesic with large winding high in the cusp, the winding will grow very quickly against the period, which is measured by hyperbolic length. 
The following results support this heuristic. By using word length instead of arc length (therefore insensitive to the cusps), Calegari \cite[Chapter 6.1]{Calegari} finds in this setting a Gaussian limiting distribution. In another variation of this problem, studying the asymptotic winding in homology of the geodesic flow, Guivarc'h and Le Jan \cite{GuivarchLeJan1993} showed that the limiting distribution is Gaussian when the surface is compact and Cauchy otherwise. In Proposition \ref{tails} we show that winding number functions on a compact hyperbolic surface grow at most linearly in $\ell_C$, ruling out a limiting Cauchy distribution.

Our main new result on the distribution of windings is the following equidistribution theorem.

\begin{theorem}\label{Thm:equidistribution}
Let $\G\in\mathcal{G}$ and consider the winding number given by \eqref{bullet}. Let $A\subseteq\Z$ be a set with natural density $d(A)$. Then the asymptotic density of prime geodesics with winding number in $A$ is equal to $d(A)$. Explicitly, if $\pi_A(T)$ denotes the number of prime geodesics in $\Pi(T)$ with winding number in $A$, then
\begin{align}\label{eq:geom Dirichlet}
\lim_{T\to\infty} \frac{\pi_A(T)}{\pi(T)}\ =\ d(A).
\end{align}
\end{theorem}

The analogon of this equidistribution theorem for prime geodesics in fixed homology classes was established by Petridis and Risager in \cite{PetridisRisager2008}. For winding numbers, the result (suggested to the second author by Risager) is new --- even in the case of the modular orbifold. We view (\ref{eq:geom Dirichlet}) as a geometric analogue of Dirichlet's prime number theorem. It would be very interesting to identify other natural winding numbers for which this arithmetic behavior persists. It is also not clear whether the restriction to the family $\mathcal{G}$ is necessary; we discuss this point further in Section \ref{discu}.

\section{Discussion of Results and Structure of the Paper}\label{discu}
Let $G=\SL_2(\R)$ and let $\G<G$ be a discrete subgroup. We say that $\G$ is a cofinite Fuchsian group if the orbifold $\G\bk\h$ is noncompact and of finite area $V$. Section \ref{sec:background} reviews classical background on Fuchsian groups, automorphic forms, and their spectral theory for the convenience of the reader and to fix notation. After this, the paper is structured in three parts.

\subsection{Rademacher symbols}
Section \ref{sec:modulargeodesics} reviews key properties of the classical Rademacher symbol  as well as various topological interpretations, following in particular Atiyah \cite{Atiyah1987} and Ghys \cite{GhysICM}. The description of $\Psi(C)$ as an alternating sum encoding the number of (oriented) turns around the cusp of the modular orbifold, seen at the opening of the introduction, does not appear in the literature and we discuss it in detail in Proposition \ref{prop:CF} and the remarks that follow it. In Section \ref{Rad}, we review properties of the Rademacher symbols for Fuchsian groups $\G$ introduced in \cite{Burrin2022}. We in particular prove in Theorem \ref{prop:Psiclassinvariant} and Corollary \ref{conj} that the Rademacher symbol $\Psi$ for $\G$ is conjugacy class invariant on the set of hyperbolic elements and hence $\Psi$ can be viewed as an invariant of the set $\Pi$ of prime geodesics on $M=\G\bk\h$. For an oriented closed geodesic $C$ on $M$, the value of the Rademacher symbol can be expressed by the real period
\begin{align}\label{period}
\Psi(C) = \int_C E_2(z)dz,
\end{align}
where the closed 1-form $E_2(z)dz$ is determined by a weight 2 nonholomorphic Eisenstein series; see \cite[Lemma 3.1]{Burrin2022}. When $\G$ is distinct from $\SL_2(\Z)$ the period \eqref{period} is not necessarily an integer --- and hence not a winding number --- and can in fact be irrational. The study of these periods is related to classical results in arithmetic geometry --- the theorems of Manin--Drinfeld and Manin--Mumford --- and this is the object of \cite{Burrin2022}. For our purposes, we only record  that the period \eqref{period} is rational whenever $\G\in\mathcal{G}$; see \cite[Theorem 1.2]{Burrin2022}.

\subsection{Winding numbers}

For ease of notation we assume in this discussion that $\G$ is torsionfree.  Geometrically this means that the quotient space $M=\G\bk\h$ has the structure of a hyperbolic surface. However our results also hold for orbifolds. The main contribution of this article is the construction of a  nowhere-vanishing (non)holomorphic automorphic form $\triangle^*_2$ of weight 2 for $\G$ such that $\triangle_2^*(z)\to0$ as $y\to\infty$. We lift this form to the continuous function 
$$ 
\widetilde\triangle_2^*:\G\bk G \to\C, \quad
 \widetilde\triangle^*_2(g)=\triangle_2^*(g(i)) g'(i),
 $$
 where $g(i)$ denotes the action of $g$ by fractional linear transformation. The automorphic transformation of $\triangle_2^*$ guarantees that the lift $\widetilde\triangle_2^*$ quotients through $\G\bk G\cong {\rm T}^1 M$. By construction, the lift $\widetilde\triangle_2^*$ yields the winding number
  $$
 {\rm ind}(\widetilde\triangle_2^*(C)) = \frac{1}{2\pi i} \int_{\widetilde\triangle_2^*(C)} \frac{dz}{z}.
 $$
 The construction of a nowhere-vanishing holomorphic automorphic form of arbitrary weight that vanishes at infinity can be achieved via the Kronecker first limit formula. This is well understood but not quite what we need: the resulting form typically comes twisted with a multiplier system, which stands as an obstruction to lifting the form to $\G\bk G$. We introduce an algebraic procedure, based on group cohomology and some differential topology, to perturb this form into one with trivial multiplier system, at the cost that the resulting form is no longer holomorphic.

 \begin{theorem}\label{Thm1b}
 Let $M$ be a cusped hyperbolic surface with finite area $V$. There exists a nonholomorphic  automorphic form $\triangle^*_{2}$ of weight $2$ for $\G$ that is nowhere-vanishing and for which $\triangle^*_2(z)\to0$ as $y\to\infty$.  Moreover, for a fixed generating set $\mathcal{S}=\{\g_i\}_{i\in I}$ for $\G$, there exists a conjugacy-class invariant function $\zeta:\G\to\Z$ such that 
$
\zeta(C) = {\rm ind}(\widetilde\triangle_{2}^*(C))
$
for each $C\in\Pi$, and
$$
\zeta(\g) =\zeta(\g_{i_1}\cdots \g_{i_k}) = \frac{V}{2\pi}\op{\Psi(\g) - \Psi(\g_{i_1}) - \dots - \Psi(\g_{i_k})}
$$
for each $\g\in\G$.
\end{theorem}

See Theorem \ref{thm1b} for the full statement for hyperbolic orbifolds. When $\Psi$ is rational-valued, it is possible to replace the nonholomorphic form $\triangle_2^*$ by a holomorphic form $\triangle_k$ of weight $k$ and degree $\tfrac{kV}{4\pi}$, at the expense of a nonexplicit and possibly higher weight $k$.  For $\PSL_2(\Z)$, $\Psi = {\rm ind}(\widetilde \triangle_{12})$, where $\triangle_{12}$ is precisely the modular discriminant.

\begin{theorem}\label{Thm1c}
Let $M$ be a cusped hyperbolic orbifold with finite area $V$. If the Rademacher symbol $\Psi$ is rational-valued, then there exists $k\in2\N$ such that 
\begin{align}\label{good-winding}
 {\rm ind}(\widetilde\triangle_k(C)) = \frac{kV}{4\pi}\Psi(C)
\end{align}
for each $C\in\Pi$.
\end{theorem}

We recall that by our previous discussion of Rademacher symbols, this situation arises for example when $\G$ is genus 0, maximal arithmetic, or a congruence subgroup of the modular group. In particular the assumption in Theorem \ref{Thm1c} is satisfied for all $\G\in\mathcal{G}$. The fact that the winding number \eqref{good-winding} is defined relative to a {\em holomorphic} form is of particular importance. It provides access to the spectral theory of automorphic forms; see Proposition \ref{prop:base ef}. 

\subsection{Counting results}
To prove Theorem \ref{Thm1} we first need the following twisted prime geodesic theorem for multiplier systems of arbitrary real weight, which we prove in Section \ref{sec:PGT}.

\begin{theorem}\label{Thm2}
Let $\G$ be a cofinite Fuchsian group and let $\chi$ be a multiplier system of weight $r$ on $\G$. Then as $T\to\infty$, we have
\begin{align}\label{tpgt-intro}
\sum_{C_\g\in \Pi(T)} \chi(\g) =  \Li\left(e^{s_0(\chi,r)T}\right)+\dots+\Li\left(e^{s_k(\chi,r)T}\right) + O(e^{3T/4}L(\chi,r)),
\end{align}
where $s_0(\chi,r)\geq s_1(\chi,r) \geq \dots\geq s_{k}(\chi,r)>\tfrac12$ are the spectral eigenparameters determined by the small eigenvalues $\lambda_j(\chi,r)=s_j(\chi,r)(1-s_j(\chi,r))$ of the weight $r$ Laplacian $\Delta_r$, and $L(\chi,r)$ is defined by Equation \ref{L}.
\end{theorem}
For $r=0$ and a multiplier system $\chi$ that is trivial at $r=0$, we recover the usual prime geodesic theorem; see, e.g., \cite[Theorem 1.5]{Iwaniec}. The proof of Theorem \ref{Thm2} follows a routine application of Selberg's trace formula for arbitrary real weight as treated by Hejhal \cite{HejhalVol2}. To keep track of the dependence of the error term on $(\chi,r)$ we derive in Section \ref{sec:Weyl} upper bound estimates for the two spectral terms appearing in Weyl's law, which may be of independent interest; see Theorem \ref{thmFvE}.

\vspace{.4cm}

With Theorem \ref{Thm2} in hand, we now review the strategy of Sarnak \cite{Sarnak2010} adapted to our setting. By construction, the winding number $\zeta(C)$ in Theorem \ref{Thm1b} defines a continuous family $(\chi_r)_{r\in\R}$ of multiplier systems via
$
\chi_r(\g) = e^{i\pi r\,\zeta(\g)}.
$
Fix $n\in\Z$. To pick up all prime geodesics $C\in\Pi(T)$ with $\zeta(C)=n$ windings, we will want to integrate (\ref{tpgt-intro}) as follows: 
\begin{align}\label{integration}
\pi_n(T) = \sum_{C\in\Pi(T)}  \int_{-1/2}^{1/2} e^{i\pi r\op{\zeta(C)-n}}\, dr = \int_{-1}^1 \Li\op{e^{s_0(\chi,r)T}} e^{-i\pi r n}\, dr + \dots.
\end{align}
There are two obstacles to carry out this strategy in our setting. The first obstacle is that we don't {\em a priori} understand how the eigenparameters $s_j(\chi,r)$ vary as functions of $r$. In Section \ref{sec:stab} we prove the existence of a small interval $|r|\leq \delta$ in which the bottom eigenparameter $s_0(\chi,r)$ has multiplicity 1; see Corollary\ref{conclu}. Observe that the width of this interval is reflected in the growth exponent of the error term in Theorem \ref{Thm1}. This result is obtained by studying the counting function for small eigenvalues of $\Delta_r$ under continuous deformation in $r$. Our strategy is inspired by the analogue for characters established by Risager \cite{Risager2011} with the addition of a crucial estimate of Jorgenson and Lundelius on the hyperbolic heat trace  \cite{JorgensonLundelius1997}. Whereas for $\PSL_2(\Z)$, one can rely on a result of Bruggeman asserting that the spectrum of $\Delta_r$ has no nontrivial small eigenvalues \cite{Bruggeman1986}, the situation for general cofinite Fuchsian groups requires the full force of the results derived in Section \ref{sec:stab}.

The second obstacle is the lack of information on the eigenparameter $s_0(\chi,r)$ to carry out the integration of the main term on the right hand side of \eqref{integration}. This issue can be circumvented if we restrict to surfaces for which $\G\in\mathcal{G}$ and work with the winding number (\ref{good-winding}) relative to the holomorphic form $\triangle_k$. Then it follows from Proposition \ref{prop:base ef} that  $s_0(\chi_r,r)= 1-\tfrac{|r|}{2}$, and Theorem \ref{Thm1} follows. 

\vspace{.2cm}
The proofs of Theorem \ref{Thm1}, Theorem \ref{Thm:Cauchy} and  Theorem \ref{Thm:equidistribution} are contained in Section \ref{sec:stat}. It would be interesting to approach these statements through different tools, in particular with the aim of extending them to the general family of winding numbers of Theorem \ref{Thm1b}.

\section{Background}\label{sec:background}

\subsection{Hyperbolic geometry}\label{2.1}
We will work on the upper half plane $\h=\{z=x+iy\in\C: y>0\}$. For each $z\in\h$, we have ${\rm T}_z\h\cong\C$ and the hyperbolic metric is given by the inner product $\scal{v,w}_z=y^{-2} \re(v\overline{w})$. 
Geodesics in this model are either vertical half-lines or semi-circles orthogonal to the real axis. The unit tangent bundle is given by ${\rm T}^1\h = \{(z,v)\in\h\times\C: |v|=y\}$. We parametrize elements of ${\rm T}^1\h$ by $(z,\theta)$ where $\theta\in\R/2\pi\Z$ is the angular variable of $v$ as measured counterclockwise from the vertical. Each point $(z,\theta)\in {\rm T}^1\h$ determines a unique hyperbolic geodesic on ${\rm T}^1 \h$ by parallel transport. The unit tangent bundle ${\rm T}^1\h$ may further be algebraically identified with the matrix group $\PSL_2(\R)$ via the transitive action on ${\rm T}^1\h$ given by
$$
(z,\theta) \mapsto (g(z), \theta - 2\arg j(g,z)),
$$
where 
\begin{align}\label{1}
g(z) &=\frac{az+b}{cz+d};\nonumber\\  j(g,z) &=cz+d;\\ 
{\rm arg}(z) &\in(-\pi,\pi].\nonumber
\end{align}
The function $j(\cdot,z)$ defines a multiplicative cocycle on $\SL_2(\R)$; i.e., it satisfies $j(gh,z)=j(g,hz)j(h,z)$ for all $g,h\in \SL_2(\R)$. Since $\det(g)=1$ and $z\in\h$, it is moreover nowhere vanishing. To see that the action is well defined, we remark that for all $g,h\in \SL_2(\R)$, the difference
\begin{align}\label{omega}
\omega(g,h) \colonequals \frac{1}{2\pi}\op{\arg j(g,hz)+\arg j(h,z)-\arg j(gh,z)}
\end{align}
is an integer in $\{0,\pm 1\}$; indeed since $\omega(g,h)$ is continuous as a function of $z$ and integer-valued, it does not depend on the particular choice of $z$. The choice of the branch of logarithm $\arg(z)\in(-\pi,\pi]$ is recorded for instance by the fact that $\omega(-1,-1)=1$. For further reference, we also record that $\omega$ defines a 2-cocycle --- indeed, direct computation establishes that
$$
\omega(gh,k)+\omega(g,h) = \omega(g,hk)+\omega(h,k)
$$
for all $g,h,k\in \SL_2(\R)$ --- and indeed is a bounded 2-cocycle representative of the Euler class \cite{Burrin2018}.

\subsection{The geometry of Fuchsian groups}
By uniformization, every hyperbolic surface (or orbifold) may be realized as the space of orbits $\G\bk\h$ for some discrete subgroup $\G<\PSL_2(\R)$ (under the action of $\G$ by fractional linear transformations as in (\ref{1})). This action has discrete orbits, which allows to represent the orbit space $\G\bk\h$ visually by a choice of fundamental domain $\cF$ in $\h$. The hyperbolic area ${\rm vol}(\cF) = \int_\cF d\mu(z) = \int_{\cF} \tfrac{dxdy}{y^2}$ does not depend on the particular choice of the fundamental domain and we set
$$
V={\rm vol}(\cF).
$$ 
We will say that $\G$ is a {\bf cofinite Fuchsian group} if $\G<\PSL_2(\R)$ is a discrete subgroup that admits a noncompact but finite-area fundamental domain in $\h$.

Elements of $\G$ are classified by their fixed point sets in $\h$; each $\g$ is either the identity, elliptic (with one fixed point in $\h$), parabolic (with one fixed point in $\partial_\infty\h=\R\cup\{\infty\})$, or hyperbolic (with two distinct fixed points in $\partial\h_\infty)$, and their conjugacy classes correspond to, respectively, the identity, torsion points of $\G\bk\h$, cusps of $\G\bk\h$, and closed  geodesics on $\G\bk\h$.  A Fuchsian group is cofinite if and only if it has finite area $V$ and a finite (but nontrivial) number of ($\G$-inequivalent) cusps. It can be presented as a finitely generated group generated by hyperbolic motions $a_1,b_1,a_2,b_2,\dots,a_g,b_g$ (with $g$ the genus of $\G\bk\h$), parabolic motions $c_1,\dots,c_h$ (with $h$ the number of cusps of $\G\bk\h$), and elliptic motions $e_1,\dots, e_\ell$ of orders $m_i\geq2$, $i=1,\dots,\ell$, satisfying the relations
$$
[a_1,b_1][a_2,b_2]\cdots[a_g,b_g] c_1\cdots c_h e_1\cdots e_\ell = e_i^{m_i} =1
$$
for each $i=1,\dots,\ell$. Each cofinite Fuchsian group is in particular the free product of a free group and finite cyclic groups. Via the Gauss--Bonnet formula we have the identity
$$
V = 2\pi\op{2g-2+h+\sum_{i=1}^\ell \op{1-m_i^{-1}}}.
$$

For arithmetic applications, it is convenient to work with discrete subgroups $\G<\SL_2(\R)$ instead. The action of $\G$ on $\h$ factors through its image in $\PSL_2(\R)$, which is given by either $\G$ or $\G/\{\pm I\}$ depending on whether $-I\in \G$. This distinction requires some caution. For instance, if $-I\in\G$, we only have a bijection between the set of closed oriented geodesics and the conjugacy classes of hyperbolic elements $\g\in\G$ with $\tr(\g)>2$. More explicitly, let $\g=\bsm a&b\\c&d\esm$ be a primitive hyperbolic matrix in $\SL_2(\R)$; then $(a+d)^2>4$ and $c\neq0$. The fixed points of $\g$ in $\partial_\infty\h$ are given by 
$$
\alpha =\frac{a-d+\sqrt{(a+d)^2-4}}{2c},\quad \overline\alpha =\frac{a-d-\sqrt{(a+d)^2-4}}{2c}.
$$
As a  hyperbolic motion, $\g$ acts on the unique hyperbolic geodesic connecting $\alpha$ and $\overline\alpha$ via the geodesic flow, and the orientation of the resulting geodesic depends on $\sign(c(a+d))$. In fact, a computation of eigenvalues and eigenvectors for $\g$ yields the diagonal form
$$
a_\g = g^{-1}\g g= \bpm \lambda&0\\ 0&\lambda^{-1}\epm,
$$
where 
$$
g=\bpm \alpha &\overline\alpha\\ 1&1\epm \quad \text{ and }\quad \lambda=\frac{a+d+\sqrt{(a+d)^2-4}}{2}.
$$
We now observe that $\lambda>1$ iff $a+d>2$ and that $\alpha>\overline\alpha$ iff $c>0$. We summarize these observations in the following table.

\begin{table}[htp]
\begin{center}
\begin{tabular}{|c|c|c|}
\hline
& $a+d>2$ & $a+d<-2$\\
\hline
$c>0$ & 
\begin{tikzpicture}
%\draw [domain={sin(2*pi/3 r)}:3.4] plot ({cos(2*pi/3 r)},\x);
\draw [domain=-2:2] plot (\x,0);
\draw  [domain=0:pi] plot ({cos(\x r)},{sin(\x r)});
\draw[] ({1)},{0}) node[below] {$\alpha$} ;
\draw[] ({-1},{0}) node[below] {$\overline\alpha$} ;
\draw[] (0,1) node[above] {$\overset{\g}{\rightarrow}$};
\end{tikzpicture}
&
\begin{tikzpicture} 
\draw [domain=-2:2] plot (\x,0);
\draw  [domain=0:pi] plot ({cos(\x r)},{sin(\x r)});
\draw[] ({1)},{0}) node[below] {$\alpha$} ;
\draw[] ({-1},{0}) node[below] {$\overline\alpha$} ;
\draw[] (0,1) node[above] {$\overset{\g}{\leftarrow}$};
\end{tikzpicture}
 \\
 \hline
$c<0$ 
&
\begin{tikzpicture} 
\draw [domain=-2:2] plot (\x,0);
\draw  [domain=0:pi] plot ({cos(\x r)},{sin(\x r)});
\draw[] ({1)},{0}) node[below] {$\overline\alpha$} ;
\draw[] ({-1},{0}) node[below] {$\alpha$} ;
\draw[] (0,1) node[above] {$\overset{\g}{\leftarrow}$};
\end{tikzpicture}
&
\begin{tikzpicture} 
\draw [domain=-2:2] plot (\x,0);
\draw  [domain=0:pi] plot ({cos(\x r)},{sin(\x r)});
\draw[] ({1)},{0}) node[below] {$\overline\alpha$} ;
\draw[] ({-1},{0}) node[below] {$\alpha$} ;
\draw[] (0,1) node[above] {$\overset{\g}{\rightarrow}$};
\end{tikzpicture}
\\
\hline
\end{tabular}
\end{center}
\captionof{table}{}\label{fig2}
\end{table}

The subgroup $\G<\SL_2(\R)$ is discrete if and only if its image in $\PSL_2(\R)$ is discrete and we keep the terminology of cofinite Fuchsian group here as well.

\subsection{Multiplier systems}\label{sec:ms}

The background material in this section follows \cite{Roelcke1966} and \cite{HejhalVol2}.
\begin{definition}
Let $r\in\R$. A function $\chi:\Gamma\to\C$ is a {\em multiplier system of weight $r$} if it verifies
\begin{itemize}
\item[(i)] $|\chi(\gamma)|=1$;
\item[(ii)] $\chi(-I)=e^{-\pi ir}$ (if $-I\in\G$);
\item[(iii)] and $\chi(\gamma_1\gamma_2)\chi(\gamma_1)^{-1}\chi(\gamma_2)^{-1}=e^{2\pi ir\omega(\gamma_1,\gamma_2)}$.
\end{itemize}
\end{definition}

If $r=n\in\Z$ is even (respectively odd), then $\chi$ is an even (respectively odd) unitary character of $\Gamma$. By periodicity, any multiplier system $\chi$ is of weight $r$ if and only if it is of weight $r+2m$ for any $m\in\Z$. Thus up to multiplication by a character, we may assume that a multiplier system has weight $r\in(-1,1]$. Further, given a multiplier system $\chi$ of weight $r$ on $\Gamma$ and an element $g\in \GL_2(\R)$, 
$$
\chi^g(\tau) \colonequals \chi(g\tau g^{-1}) e^{2\pi ir(\omega(g\tau g^{-1},g)-\omega(g,\tau))}
$$
defines a multiplier system on $g^{-1}\Gamma g$. 

\subsection{Spectral theory}\label{sec:specth}

Let $\chi$ be a multiplier system of weight $r$ for $\Gamma$, and let $\mathcal{H}(\G,\chi,r)$ be the space of all measurable functions $f:\h\to\C$ such that 
\begin{align}\label{real wt}
f(\gamma z) =\chi(\gamma)\left(\frac{j(\gamma,z)}{|j(\gamma,z)|}\right)^{r} f(z)
\end{align}
and
\begin{align*}
\scal{f,f} &= \int_{\Gamma\backslash\h} |f(z)|^2 dz < +\infty. 
\end{align*}
Equipped with the inner product $\scal{\cdot,\cdot}$, $\cH(\G,\chi,r)$ has the structure of a Hilbert space and is isomorphic to the space $L^2(\Gamma\backslash\h,\chi,r)$ of (equivalence classes of) square-integrable functions satisfying \eqref{real wt}. Acting on (a dense subspace of) $L^2(\G\bk\h,\chi,r)$ we have the Maass raising and lowering operators
\begin{align*}
K_r &= iy\partial_x + y\partial_y + \frac{r}{2},\\
\Lambda_r &= iy \partial_x - y\partial_y+\frac{r}{2}.
\end{align*}
For $f\in\mathcal{H}(\G,\chi,r)$ and $g\in \mathcal{H}(\G,\chi,r+2)$ both at least $C^1$, we have the identity
$$
\scal{K_r f,g} = \scal{f,\Lambda_{r+2}g}.
$$
The Maass operators relate to the weight $r$ (geometric) Laplacian via
\begin{align*}
\Delta_r &= \Lambda_{r+2} K_r -\frac{r}{2}\left(1+\frac{r}{2}\right) = K_{r-2}\Lambda_r +\frac{r}{2}\left(1-\frac{r}{2}\right) =  -y^2 \left(\partial_x^2 +\partial_y^2\right) +iry\partial_x.
\end{align*}
It follows that 
\begin{align}\label{31aout}
\scal{f,\Delta_r g} = \scal{K_r f, K_r g} -\frac{r}{2}\left(1+\frac{r}{2}\right)\scal{f,g} = \scal{\Lambda_r f,\Lambda_r g}+\frac{r}{2}\left(1-\frac{r}{2}\right)\scal{f,g}
\end{align}
and that 
$$
\scal{\Delta_r f,g}=\scal{f,\Delta_r g}
$$ 
for all $C^2$ functions $f,g\in \cH(\G,\chi,r)$.  

If $-I\not\in\G$, then each cusp $\fa$ has an infinite cyclic stabilizer group $\G_\fa<\G$. We will denote its generator by $\gamma_\frak{a}$. If $-I\in \G$, then $\G_\fa$ is isomorphic to $\Z\times\Z/2\Z$; in this case we will also denote the infinite cyclic generator by $\g_\fa$. If $\fa,\fb$ are equivalent cusps, then $\G_\fa$ and $\G_\fb$ are conjugate in $\G$. Let $\fa$ be a cusp for $\G$ and let $\g_\fa$ be a cyclic generator of $\G_\fa$. We say that $\fa$ is {\bf singular} if $\chi(\g_\fa)=1$ and {\bf regular} otherwise. With respect to (the unique self-adjoint extension of) $\Delta_r$, the space $L^2(\Gamma\backslash\h,\chi,r)$ has a complete spectral resolution with pure point spectrum
$$
\lambda_0(\chi,r)\leq \lambda_1(\chi,r)\leq \lambda_2(\chi,r) \leq \dots\quad \to\infty
$$
and absolutely continuous spectrum $[\tfrac14,\infty)$ with multiplicity $m=m(\G,\chi,r)$ equal to the number of inequivalent singular cusps. The continuous spectrum is described by Eisenstein series. 

The bottom eigenvalue is given by
$$
\lambda_0(\chi,r) = \inf_{f\in\cH\cap C^2} \frac{\scal{\Delta_r f,f}}{\scal{f,f}}.
$$
Specializing \eqref{31aout} to $f=g$ yields
$$
\scal{\Delta_r f,f} = \| K_r f\|^2 -\frac{r}{2}\left(1+\frac{r}{2}\right)\|f\|^2 = \|\Lambda_r f\|^2 +\frac{r}{2}\left(1-\frac{r}{2}\right)\|f\|^2,
$$
and implies that
$$
\|K_r f\|^2 = \|\Lambda_r f\|^2 +r\|f\|^2.
$$
That is, $\|K_r f\|^2 \geq \|\Lambda_r f\|^2$ if $r\geq 0$ and $\|\Lambda_r f\|^2\geq \|K_r f\|^2$ if $r\leq0$. Hence 
$$
\lambda_0(\chi,r) = \frac{|r|}{2}\op{1-\frac{|r|}{2}} + \inf_{f\in\cH\cap C^2} \begin{dcases} \frac{\|\Lambda_r f\|^2}{\|f\|^2} & \text{ if } r\geq0,\\ \frac{\|K_r f\|^2}{\|f\|^2} &\text{ if }r\leq 0.
\end{dcases}
$$

In Section \ref{Rad}, we will exhibit a continuous family $\{\chi_r\}_{r\in\R}$ of multiplier systems (for any cofinite Fuchsian group) and show that for each multiplier system $\chi_r$, we have
$$
\lambda_0(\chi_r,r) = \frac{|r|}{2}\op{1-\frac{|r|}{2}},
$$
or in other words that the minimal bottom eigenvalue for $\Delta_r$ is realized.

\subsection{Eisenstein series and residual spectrum}

Computations as well as geometric considerations are easier when the cusp is at $\infty$, which we may assume up to conjugating $\G$ by a {\bf scaling} $\gs_\fa\in G$ such that $\gs_\fa(\infty)=\fa$ and $(\gs_\fa^{-1}\G\gs_\fa)_\infty=\gs_\fa^{-1}\G_\fa\gs_\fa$ is generated by 
$$
T = \bpm 1 &1\\0&1\epm
$$
if $-I\not\in\G$ or is otherwise isomorphic to $\Z\times\Z/2\Z$ with the cyclic generator of $\Z$ given by $T$. The choice of scaling is not unique; indeed, each element of the continuous family $\sigma_{\frak{a}} T^x=\gs_\fa\bsm 1&x\\0&1\esm$, $x\in\R$, showcases those same properties.

For each singular cusp $\fa$,  we have a well defined Eisenstein series given by
$$
E_\fa(z,s,\chi,r) = \sum_{\g\in\G_\fa\bk\G} \overline{\chi(\g)} e(-r\omega(\gs_\fa,\g)) e^{-ir\arg j(\gs_\fa \g,z)} \im(\gs_\fa \g z)^s,
$$
which converges absolutely and uniformly on compact sets for $\re(s)>1$, transforms according to (\ref{real wt}), verifies $(\Delta_r-s(s-1)) E_\fa(z,s,\chi,r)=0$ and has a meromorphic continuation to all $s\in\C$.

The $m\times 1$ vector $\mathcal{E}(z,s,\chi,r)$ of all Eisenstein series at (inequivalent) singular cusps satisfies the functional equation
\begin{align}\label{functionaleq}
\mathcal{E}(z,1-s,\chi,r) = \Phi(s,\chi,r)\mathcal{E}(z,s,\chi,r),
\end{align}
where $\Phi=\left(\varphi_{\frak{ab}}\right)_{\fa,\fb}$ is the {\bf scattering matrix}, and $\varphi=\det\Phi$ is the {\bf scattering determinant}. From the functional equation (\ref{functionaleq}), we see that 
$$
\Phi(s,\chi,r)\Phi(1-s,\chi,r)=I.
$$ 
Furthermore, on the line $\re(s)=\tfrac12$, we have
$$
\Phi\left(\frac12+it,\chi,r\right)\Phi\left(\frac12+it,\chi,r\right)^*\ =\ \Phi\left(\frac12+it,\chi,r\right)\Phi\left(\frac12-it,\chi,r\right)=I,
$$
where $A^*=\overline{A}^T$ denotes the conjugate transpose of the matrix $A$. In particular, the entries of the scattering matrix are bounded. Following Selberg \cite[pp.~655-656]{Selberg1989}, the scattering determinant $\varphi=\det\Phi$ has  the Dirichlet series expression
\begin{align}\label{Dirich}
\varphi(s,\chi,r) = \left(\frac{\sqrt\pi 4^{1-s}\Gamma(2s-1)}{\Gamma(s+\tfrac{r}{2})\Gamma(s-\tfrac{r}{2})}\right)^{m(\chi,r)}\sum_{n\geq 1}\frac{a_n}{b_n^{2s}}
\end{align}
where $0<b_1<b_2<\dots$ are real coefficients, and where the series converges absolutely for $\re(s)>1$, and is moreover regular for $\re(s)\geq1/2$ except possibly for a finite number of poles $s_1(\chi,r), s_2(\chi,r),\dots, s_k(\chi,r)$ in the interval $(1/2,1]$, each with multiplicity bounded above by $m(\chi,r)$. These poles are the eigenparameters of the residual eigenvalues 
$$
\lambda_j(\chi,r)=s_j(\chi,r)(1-s_j(\chi,r))
$$ 
of the weight $r$ Laplacian $\Delta_r$. 

If $E_\frak{a}(z,s,\chi,r)$ has a pole at $s=s_0$, this pole is simple, $s_0\in(1/2,1]$, and $\varphi_\frak{aa}(s)$ also has a pole at $s_0$, see \cite[S\"atze 10.3-4]{Roelcke1966}. Conversely, if $\varphi_\frak{aa}(s)$ has a pole at $s_0\in(1/2,1]$, then the Maass--Selberg relations (see \cite[Lemma 11.2]{Roelcke1966}) imply that $E_\frak{a}(z,s,\chi,r)$ has a pole at $s_0$. The residue $\rho_{\frak{a},s_0}\in L^2(\Gamma\backslash\h,\chi,r)$ at this pole is an eigenfunction of $\Delta_r$ with eigenvalue $s_0(1-s_0)$. Moreover, if $\rho_{\frak{b},s_0}$ is the residue of $E_\frak{b}(z,s,\chi,r)$ at $s=s_0$, then 
\begin{align}\label{2.26}
\scal{\rho_{\frak{a},s_0},\rho_{\frak{b},s_0}} = \operatorname{Res} \varphi_\frak{ab}(s)\vert_{s=s_0},
\end{align}
 see \cite[Satz 11, p.302]{Roelcke1966}. Finally, we record the following lemma.

\begin{lemma}\label{lm:mult}
If the scattering determinant $\varphi$ has a pole of order $n$ at $s=s_0$, where $\re(s_0)>\tfrac12$, then there are $n$ linearly independent eigenfunctions of $\Delta_r$ with eigenvalue $s_0(1-s_0)$. 
\end{lemma}

\begin{proof} 
Suppose that for some subset $\fa_1,\dots,\fa_n$ of singular cusps, the residues $\rho_{\fa_i,s_0}(z)$ are not identically zero. By \eqref{2.26}, we can transform $\Phi(s,\chi,r)$ via elementary row operations into a matrix that only has poles at $s=s_0$ in rows $\fa_1,\dots,\fa_n$. Then $\varphi(s,\chi,r)$ has a pole of order at most $n$ at $s=s_0$. Since the order is $n$, the residues $\rho_{\fa_i,s_0}$ are linearly independent with eigenvalue $s_0(1-s_0)$.
\end{proof}

\subsection{Selberg's trace formula}\label{subsec:STF}

We consider again $\chi$ to be a multiplier system of weight $r$. We review Selberg's trace formula for $L^2(\Gamma\backslash\h,\chi,r)$ for the convenience of the reader following \cite[pp.~412-413]{HejhalVol2}. 

Let $h:\C\to\C$ be a holomorphic even function defined on 
$$
\left\{z\in\C : |y|<\max\{\tfrac{|r|-1}{2},\tfrac12\}+\delta\right\}
$$
for some $\delta>0$ such that $h(t)=O\left((1+t)^{-2-\delta}\right)$ and let 
$$
g(t) = \frac{1}{2\pi} \int_{-\infty}^\infty h(u)e^{-itu}du.
$$
Then
\begin{align}
\sum_{n\geq0} h(t_n) = &
\frac{V}{4\pi }\int_{-\infty}^\infty th(t)\frac{\sinh(2\pi t)}{\cosh(2\pi t)+\cos(\pi r)} dt\label{II}
\\
& + \frac{V}{4\pi}\sum_{\substack{l \text{ odd}\\ 1\leq l\leq|r|}} (|r|-l)h\left(\frac{i(|r|-l)}{2}\right)\label{I}
\\
& + \sum_{\substack{\{\gamma\}\\ \tr(\gamma)<2\\ 0<\theta_\gamma<\pi}} \frac{\chi(\gamma)ie^{i(r-1)\theta}}{4|\gamma| \sin\theta_\gamma} \int_{-\infty}^\infty g(t)e^{(r-1)t/2} \frac{(e^t-e^{2i\theta})}{\cosh t -\cos 2\theta_\gamma}dt\label{E}
\\
& + 
\sum_{\substack{\{\gamma\}_{\rm pr}\\ \tr(\gamma)>2}} \sum_{k\geq1} \frac{\chi(\gamma^k)\ell_\gamma g(k\ell_\gamma)}{\sinh(k\ell_\gamma/2)}\label{H}
\\
&  -g(0)\sum_{\alpha_\frak{a}(\chi,r)\neq0} \log|1-e(\alpha_\frak{a}(\chi,r))|\label{A}\\
& + \frac12\sum_{\alpha_\frak{a}(\chi,r)\neq0} (\tfrac12-\alpha_\frak{a}(\chi,r)) \mathrm{PV}\int_{-\infty}^\infty g(t)e^{(r-1)t/2} \frac{(e^t-1)}{\cosh(t)-1} dt \label{B}
\\
& + 
 m(\chi,r)\int_0^\infty \frac{g(t)(1-\cosh(rt/2))}{e^{t/2}-e^{-t/2}}dt \label{C}\\
 & -m(\chi,r)\left( g(0)\log 2 -\frac{1}{2\pi} \int_{-\infty}^\infty h(t)\frac{\Gamma'(1+it)}{\Gamma(1+it)}dt\right)\label{D}\\
 & + \frac{1}{4}h(0)\mathrm{Tr}\left(I-\Phi(\tfrac12,\chi,r)\right) + \frac{1}{4\pi} \int_{-\infty}^\infty h(t)\frac{\varphi'}{\varphi}\left(\frac12+it,\chi,r\right)dt,\label{S}
 \end{align}
 where 
 \begin{itemize}
 \item the sum on the LHS runs over the spectral eigenparameters $t_n=t_n(\chi,r)\in\C$ given by
 \begin{align}\label{param}
\lambda_{n}(\chi,r) &= s_{n}(\chi,r)(1-s_{n}(\chi,r))\\
&= \left(\frac12+it_{n}(\chi,r)\right)\left(\frac12-it_{n}(\chi,r)\right) = \frac14 + t_{n}(\chi,r)^2
\end{align}
and  counted with multiplicity;
 \item $|\g|$ in (\ref{E}) is the order of the elliptic element $\g$, and $\theta_\g$ is  the angle in $(0,2\pi)$ for which $\bsm \cos\theta_\g&-\sin\theta_\g\\ \sin\theta_\g&\cos\theta_\g\esm$ is a $\SL_2(\R)$-conjugate of $\gamma$;
 \item $\ell_\gamma$ in (\ref{H}) denotes the hyperbolic length of the unique closed geodesic corresponding to the conjugacy class $\{\gamma\}_{\rm pr}$ of the primitive hyperbolic element $\g$;
 \item the sums in (\ref{A}), (\ref{B}) are taken over the set of regular cusps with $\alpha_\fa(\chi,r)\in[0,1)$ determined by $\chi(\g_\fa)=e(\alpha_\fa)$;
 \item and the terms (\ref{C}), (\ref{D}), (\ref{S}) only appear in the presence of singular cusps, whereby $m(\chi,r)$ is the number of inequivalent singular cusps.
 \end{itemize}

\section{Classical Rademacher Symbols and Windings of Modular Geodesics}\label{sec:modulargeodesics}

\subsection{Dedekind sums and the classical Rademacher symbol}\label{subsec:DedRad}

In this section we fix $\G=\PSL_2(\Z)$. In the classical theory of modular forms, the modular discriminant function $\triangle(z)$ is the prototype of a holomorphic cusp form for $\G$. Jacobi proved that $\triangle$ admits the infinite product expansion
$$
\triangle(z) = q\prod_{n\geq1}(1-q^n)^{24}
$$
(with $q=e^{2\pi iz}$ and $\im(z)>0$). 
It follows that $\triangle$ is nowhere-vanishing with $\triangle(z)\to0$ as $y\to\infty$.  We fix the branch of logarithm
$$
\log\triangle(z) = 2\pi iz -24\sum_{m,n\geq1}\frac{q^{mn}}{n}.
$$
Then a simple direct computation shows that
\begin{align}\label{diffeq}
\frac{d}{dz} \log\triangle(z) &=
2\pi i\op{ 1 -24 \sum_{n\geq1} \op{ \sum_{d\mid n} d} q^n} =
2\pi i\op{ E_2(z) +\frac{3}{\pi y}},
\end{align}
where $E_2$ is Hecke's modular, nonholomorphic, Eisenstein series of weight two, given by
$$
E_2(z) = \lim_{\eps\to0^+} \sum_{\g\in\G_\infty\bk\G} \frac{1}{(cz+d)^2|cz+d|^\eps}.
$$
For each $\g=\bsm a& b\\ c&d\esm\in\PSL_2(\Z)$, we have the transformation law
\begin{align}\label{triangolo}
\log\triangle(\gamma z) - \log \triangle(z)\ =\ 6\log(-(cz+d)^2) +2\pi i\Phi(\gamma),
\end{align}
where, on the right-hand side (RHS), $\log$ denotes the principal branch of logarithm and the {\bf Dedekind symbol} $\Phi(\gamma)$ is explicitly determined by the matrix entries of $\gamma$ and involves Dedekind sums \cite{Dedekind1892}. The transformation law can be deduced either from \eqref{diffeq} or via Kronecker's first limit formula; see \cite[Chapter 1.2]{Siegel}. In his study of Dedekind sums, Rademacher \cite{Rademacher1956} put the emphasis on a second invariant, the {\bf Rademacher symbol} $\Psi$,  given by 
\begin{align*}
\Psi(\gamma)\ =\ \Phi(\gamma)-3\ \sign(c(a+d)),
\end{align*}
which he shows to be integer-valued and conjugacy class invariant. (We recall that $\sign(c)=\tfrac{c}{|c|}$ if $c\neq0$ and $\sign(c)=0$ if $c=0$.)

\subsection{Topological realizations}
Several topological interpretations of $\Psi$ were identified by Atiyah \cite[(5.60)]{Atiyah1987}. For instance, for $\gamma\in \SL_2(\Z)$ hyperbolic, the normalization $-\tfrac13 \Psi(\gamma)$ coincides with Meyer's function in the study of signatures of torus bundles over surfaces, with Hirzebruch's signature defect, with the logarithmic monodromy of Quillen's determinant line bundle, or with the Atiyah--Patodi--Singer eta-invariant %(unrelated to Dedekind's $\eta$) 
and its adiabatic limit, to name just a few. 

Closer to our interest in this paper, $\Psi(\gamma)$ also coincides (again for $\gamma$ hyperbolic, and possibly up to  a scalar multiple) with the signature of an oriented link on the braid group $B_3$ with three strands \cite{GambaudoGhys2005}, with the linking number for modular knots with the trefoil knot $\tau$ in $S^3\setminus\tau$, and with the planar winding number
\begin{align}\label{planar}
{\rm ind}(\widetilde\triangle(C),0)\ =\ \frac{1}{2\pi i}\int_{C} \frac{d\widetilde\triangle}{\widetilde\triangle},
\end{align}
where $\widetilde\triangle(C)\subset \C\setminus\{0\}$ is the continuous closed curve obtained by taking the lift $\widetilde\triangle:\SL_2(\Z)\bk \SL_2(\R)\to\C^*$ of the modular discriminant $\triangle(z)$, and evaluating it along a parametrized closed oriented orbit $C$ of the geodesic flow on $\G\bk G$ \cite{GhysICM}. The equivalence of these last three interpretations boils down to the remarkable fact that the homogeneous space $\SL_2(\Z)\bk\SL_2(\R)$ is diffeomorphic to $S^3\setminus\tau$ with common fundamental group $B_3$; see \cite[p.~84]{Milnor}. 

\subsection{Rademacher symbol as winding number}
In the rest of this section, we give an explicit interpretation of the winding number $\Psi$, which corresponds to the description given in the introduction. The explicit formula for the Rademacher symbol in terms of Dedekind sums can be evaluated via the Euclidean algortihm, and thus in terms of continued fraction expansions; see \cite{Meyer1957,Zagier1975,Hickerson1977,Kelmer2012}.

Let $\g\in\SL_2(\Z)$ be a primitive hyperbolic element with positive trace. Its fixed points $x_+>x_-$ are quadratic irrationals, hence admit an eventually periodic simple continued fraction expansion. Let 
$$
x_+=[a_1,\dots,a_m,\overline{a_{m+1},\dots,a_{m+n}}],
$$
where $a_i\in\Z$ and $a_i\geq1$ for $i\geq2$. 
The stabilizer of the purely periodic part $x=[\overline{a_{m+1},\dots,a_{m+n}}]$ in $\GL_2(\Z)$ is an infinite cyclic group generated by
$$
B = A_{a_{m+1}}A_{a_{m+2}}\cdots A_{a_{m+n}} = \bpm a_{m+1}& 1\\ 1& 0\epm\bpm a_{m+2}&1\\1& 0\epm \cdots \bpm a_{m+n}&1\\1&0\epm.
$$
Note that if $n$ is even, we have $B\in\SL_2(\Z)$. The Rademacher function $\Psi$ behaves nicely on such group elements..
\begin{proposition}\label{prop:CF}
Let $n\in\N$ be even, and let $a_1,\dots,a_n$ be a finite set of positive integers. Then 
$$
\Psi(A_{a_1}\cdots A_{a_n}) = a_1 -a_2 + \dots -a_n.
$$
\end{proposition}
\begin{proof}
Let $S=\bsm 0&-1\\1&0\esm$, $T=\bsm 1&1\\0&1\esm$, $U=\bsm 1&0\\1&1\esm$. 
The proof will proceed by induction via repeated applications of the following formulas 
\begin{align*}
\Psi(T^a) &=a,\qquad \Psi(U^a)=\Psi(ST^{-a}S^{-1})=-a;\\
\Phi(\g \g') &= \Phi(\g) + \Phi(\g') -3\sign(c_\g c_{\g'} c_{\g\g'});\\
\Psi(\g) &= \Phi(\g) -3\sign(c(a+d));
\end{align*}
see \cite[(59), (62)--(63)]{RademacherGrosswald}. 

We first consider $n=2$ and write down explicitly
$$
A_a A_b = \bpm ab+1 & a \\ b& 1\epm = \bpm 1 & a \\0&1 \epm \bpm 0&-1\\ 1&0\epm\bpm 1 & -b\\0&1\epm \bpm 0&1\\-1&0\epm = T^a U^b.
$$
Then 
\begin{align*}
\Psi(A_a A_b) &= \Phi(A_a A_b) - 3 = \Phi(T^a) + \Phi(U^b) -3 = \Psi(T^a) + \Psi(U^b) =  a-b.
\end{align*}
 More generally, for any even product $A_{a_1}\cdots A_{a_n}$, all entries are positive and writing $A=A_{a_3}\cdots A_{a_n}$, we have
 \begin{align*}
 \Psi(A_{a_1}\cdots A_{a_n}) &= \Psi(A_{a_1}A_{a_2}A) = \Phi(A_{a_1}A_{a_2}A)-3\\
 &=
 \Phi(A_{a_1}A_{a_2}) -3 +\Phi(A) - 3 = \Psi(A_{a_1}A_{a_2}) +\Psi(A).
 \end{align*}
 We conclude by induction.
\end{proof}

Up to replacing $\g$ by a conjugate, we may assume that the quadratic irrational $x_+$ is reduced, i.e., that $x_+>1,$ $x_-\in(-1,0)$, and as such admits a purely periodic continued fraction expansion. We review this in detail, and take the opportunity to fix some notation. Let $B=A_{a_{m+1}}\cdots A_{a_{m+n}}$ as above and set $A=A_{a_1}\cdots A_{a_m}$. We find that $\g = ABA^{-1}$. Then
\begin{itemize}
\item If both $m$ and $n$ are even, then both $A, B\in \SL_2(\Z)$ and we have that 
$$
\Psi(\g)=\Psi(B) = a_{m+1} - a_{m+2} + \dots -a_{m+n},
$$
where $(a_{m+1},\dots,a_{m+n})$ is the periodic part of the continued fraction expansion of the fixed point $x_+$ of $\g$.
\item If $n$ is even and $m$ is odd, then $\g$ and $B$ are not conjugate in $\SL_2(\Z)$. Let
$$
w= \bpm 1 &0\\0&-1\epm
$$
and define $\phi_w(A)=w^{-1}Aw$. Explicitly, we have
$$
\bpm -1 &0\\0&1\epm \bpm a& b\\ c&d\epm \bpm 1 & 0\\ 0&-1\epm = \bpm a&-b\\ -c & d\epm.
$$
Then $\g$ is conjugate to $\phi_w(B)=A_{-a_{m+1}}\cdots A_{-a_{m+n}}$ in $\SL_2(\Z)$ and hence
$$
\Psi(\g) = \Psi(\phi_w(B)) = - a_{m+1} +a_{m+2} -\dots + a_{m+n} = -\Psi(B).
$$
\item 
Considering the closed oriented geodesic $C_\g$ associated to $\g$, one sees that $\phi_w(\g)$ corresponds to $\phi(C_\g)$ where $\phi(x+iy,\theta)=(-x+iy, \pi-\theta)$.
\begin{center}
\begin{tikzpicture} 
%\draw [domain={sin(2*pi/3 r)}:3.4] plot ({cos(2*pi/3 r)},\x);
\draw [domain=-6:2] plot (\x,0);
\draw  [domain=0:pi] plot ({-4+cos(\x r)},{sin(\x r)});
\draw  [domain=0:pi] plot ({cos(\x r)},{sin(\x r)});
\draw[gray,dashed]  [domain=0:2] plot (-2,\x);
\draw[] ({-3},{0}) node[below] {-$\alpha$} ;
\draw[] ({-5)},{0}) node[below] {-$\overline\alpha$} ;
\draw[] ({1)},{0}) node[below] {$\overline\alpha$} ;
\draw[] ({-1},{0}) node[below] {$\alpha$} ;
\draw[] (-4,1) node[above] {$\overset{\phi_w(\g)}{\leftarrow}$};
\draw[] (0,1) node[above] {$\overset{\g}{\rightarrow}$};
\end{tikzpicture}
\end{center}

\item If $n$ is odd, then $B\not\in\SL_2(\Z)$ but $B^2\in\SL_2(\Z)$, and we find that $\g$ is conjugate to either $B^2$ or $\phi_w(B^2)$ (depending on the parity of $m$). As already observed by Kelmer \cite[Proposition 3.3]{Kelmer2012}, $n$ is odd exactly when $\g$ is inert, i.e.,~invariant under the orientation-reversal involution $\phi_w$. Moreover, we have that
$$
\Psi(\g) = \pm(a_{m+1}-a_{m+2} + \dots +a_{m+n} -a_{m+1} + \dots -a_{m+n})=0.
$$
\end{itemize}

There are yet other geometric interpretations. Consider $\SL_2(\Z)$ as the free product of the finite cyclic groups generated by $S$ and $TS$. Every nontrivial primitive hyperbolic element of $\SL_2(\Z)$ can be seen to be conjugate to a word of the form $A=SU^{\eps_1} S U^{\eps_2}\cdots S U^{\eps_r}$, with $\eps_j=\pm 1$. The Rademacher symbol for $A$ can be expressed as 
$$
\Psi(\g) = \sum_{j=1}^r \eps_j.
$$
See \cite[Chapter 4-C]{RademacherGrosswald} for proofs and \cite[Section B-4]{BargeGhys1992} for a geometric description of $\Psi$ recording the (oriented) turns of the corresponding geodesic in the Basse--Serre tree of $\SL_2(\Z)$. 

\section{Winding Numbers for Prime Geodesics}\label{winding}

We fix $\Gamma$ to be a cofinite Fuchsian group with standard cusp at infinity.
In a first step, we will review how to (canonically) attach to this cusp a holomorphic automorphic form satisfying the required properties that it is nowhere-vanishing on $\h$ but vanishes at infinity. This construction already appears in the literature in connection to the Kronecker first limit formula and the theory of multiplier systems \cite{Goldstein1973,Patterson1975,HejhalVol2}. We will follow the presentation of \cite{Burrin2017}. The result is almost what we want, except that the form comes twisted by a multiplier system. In a second step, we develop an algebraic approach to trivialize said multiplier system. 

We start from the standard cusp at infinity and its associated nonholomorphic Eisenstein series $E(z,s)=E_\infty(z,s,1,0)$. This is a meromorphic function on the complex plane, with Laurent expansion 
\begin{align}\label{eq:KLF}
E(z,s) = \frac{V^{-1}}{s-1}+ K(z) + O(|s-1|)
\end{align}
as $s\to1$. We can deduce from the properties of the Eisenstein series that the Kronecker limit $K:\G\bk\h\to\R$ is $\G$-invariant, real-valued, real-analytic, and almost harmonic, namely $\Delta K = -V^{-1}$. Let $f:\h\to\C$ be a holomorphic function with 
$$
\re f(z)=K(z)+V^{-1}\log(y).
$$ 
(Note that $\Delta(K(z)+V^{-1}\log y)=0$.) Fixing the branch of logarithm with $\arg(z)\in(-\pi,\pi)$, we find that for each $\g\in\G$, we have
$$
f(\g z)-f(z) = 2V^{-1} \log j(\g,z)+i S(\g),
$$
where $S:\G\to\R$ is a quasimorphism that does not depend on the particular choice of $z\in \h$ and for which
 \begin{align}\label{future}
 S(\g_1\g_2) - S(\g_1)-S(\g_2) = 4\pi V^{-1} \omega(\g_1,\g_2)
 \end{align}
 for all $\g_1,\g_2\in\G$. We refer the reader back to Section \ref{2.1} for the cohomological interpretation of this equation. Next we easily verify that for each $r\in \R$, the function
 \begin{align}\label{def:triangle}
 \triangle_r (z) \coloneqq \exp(rV f(z)/2)
\end{align}
 is holomorphic and nowhere-vanishing on $\h$, and verifies 
 $$
 \triangle_r(\g z) = \chi_r(\g) j(\g,z)^r \triangle_r(z)
 $$
 with $\chi_r(\g) = e^{ir VS(\g)}$ for all $\g\in\G$. Starting off in the Kronecker limit formula (\ref{eq:KLF}) with the Fourier expansion of Eisenstein series, we can choose $f$ so that it has $q$-expansion given by
$$
f(z) = iz + \sum_{n\geq1} c(n)q^n,
$$
where the $c(n)$ are real coefficients growing at most polynomially in $n$. These coefficients can be written down explicitly, but we will not need this information. For our purposes we only note that when $r>0$, we have $\triangle_r(z)\to0$ as $y\to\infty$. A more careful study of the Kronecker limit formula will also show that $\triangle_r$ is holomorphic at each cusp of $\G$. In summary, we have an automorphic form with the required properties for the construction of the winding number function already described in the introduction, up to the presence of the multiplier system $\chi_r$. In particular, when $\chi_r$ is nontrivial, the lift $\widetilde\triangle_r:G\to\C^*$ does not quotient through $\G$. Our goal from here on is to `trivialize' the multiplier system $\chi_r$. 

This is essentially immediate if the function $S:\G\to\R$ that determines the multiplier system is rational-valued. Note that the constant $\tfrac{V}{2\pi}$ is a rational number (e.g., by the Gauss--Bonnet theorem) and the relation (\ref{future}) implies that $\tfrac{V}{2\pi}S$ takes values in $\tfrac{1}{m}\Z$ for some minimal $m\in \N$. Then $\chi_{2m}=1$ and we think of the winding numbers ${\rm ind}(\widetilde \triangle_{2m}(C))$ as computing the winding of the closed curve $C$ about the prescribed cusp at infinity. In the case of the modular orbifold $\PSL_2(\Z)\bk\h$, we precisely recover the winding number described in the previous section. Indeed, in this case, the Fourier expansion of the Eisenstein series yields
\begin{align*}
\lim_{s\to1^+} \op{E(z,s)-\frac{3}{\pi(s-1)}} 
= c_0  - \frac{3}{\pi} \log y - \frac{1}{2\pi}\re \log\triangle(z)
\end{align*}
hence $\triangle_{12}(z)=\triangle(z)$ is the classical modular discriminant, and 
$$
{\rm ind}(\widetilde \triangle(C))=\Psi(C).
$$

This being said, we will see via Lemma \ref{lm-Psi} that the question of whether $S$ takes only rational values is answered negatively \cite{Burrin2022}. Further, a characterization of Fuchsian groups for which $S$ is rational seems out of reach. In the rest of this section, we develop an algebraic approach to trivialize the multiplier system by perturbing $\triangle_r$. The trade-off is that the modified form $\triangle^*_r$ is no longer holomorphic.

The starting point of our construction is the observation that (\ref{future}) takes values in $\Q$. Set
$$
dS(\g_1,\g_2) \coloneqq S(\g_1\g_2) -S(\g_1) -S(\g_2). 
$$ 
By definition, $dS$ is a 2-cocycle on $\G\times\G$. Since ${\rm H}^2(\G,\Q)$ is trivial, it is in fact a 2-coboundary and there exists $\zeta:\G\to\Q$ satisfying $d\zeta=dS$. Inspired by this fact, we will construct a winding number function $\zeta:\G\to\Z$ with the property that $d\zeta=dS$. We first introduce the notion of the torsion index.
\begin{definition}
We define the {\bf torsion index} of $\G$ to be the smallest $\tau\in\N$ such that $\tfrac{\tau V}{2\pi}\in\N$.
\end{definition}

The torsion index of $\SL_2(\Z)$ is $\tau=6$, while if $\G$ is torsionfree, it is $\tau=1$. Fix a generating set $\cS=\{\g_i\}_{i\in I}$ for $\G$, including $-I$ if $-I\in\G$. Set $g(\g)=-S(\g)$ for each nonelliptic generator $\g\in\cS$, $\g\neq-I$,  $g(\g)=0$ for each elliptic generator $\g\in\cS$, and $g(-I)=0$, then extend $g$ to a homomorphism on $\G$. 
We set $$\zeta\coloneqq \tfrac{\tau V}{2\pi}(S+g)$$ and claim that with this choice, $\zeta:\G\to\Z$; i.e., $\zeta$ is integer-valued. By recursion on (\ref{future}) we have
$$
S(\g^m) = m\, S(\g) +4\pi V^{-1}\sum_{k=1}^{m-1} \omega(\g^k,\g).
$$
Thus if $\g^m=\pm I$, then $S(\g)\in\tfrac{2\pi}{mV}\Z$. Then for each elliptic $\g\in\G$, we have $S(\g)\in \tfrac{2\pi}{\tau V}\Z$. It follows that each $\zeta(\g)$ is a linear combination of values of $\tfrac{\tau V}{2\pi}S$ on elliptic generators and values of $\tfrac{\tau V}{2\pi}dS$, hence an integer.

The homomorphism $\zeta - \tfrac{\tau V}{2\pi}S\in {\rm Hom}(\G,\R)$ factors through the image $\overline\G$ of $\G$ in $\PSL_2(\R)$. The classical theorems of Hurewiecz and de Rham yield the following sequence of isomorphisms:
$$
{\rm Hom}(\overline\G,\R) \cong {\rm Hom}(\overline\G/[\overline\G,\overline\G],\R) \cong {\rm Hom}({\rm H}_1(M),\R) \cong {\rm H}^1_{\rm dR}(M).$$
In particular there exists a unique smooth closed 1-form $\alpha\in\Omega^1(M)$ such that
$$
\zeta(\g) - \tfrac{\tau V}{2\pi}S(\g) = \int_{[\g]} \alpha,
$$
where on the right the integral is taken over a cycle representative of $[\g]$, the image of $\g$ under the projection $\overline\G\to \overline\G/[\overline\G,\overline\G]\cong {\rm H}_1(M)$. Let $\pi:\h\to \G\bk\h$ be the canonical covering map. On the universal cover $\h$, the pull-back form $\pi^*\alpha$ is exact (since $\h$ is contractible) and
so there exists a smooth real-valued function $a\in C^\infty(\h)$ such that 
$$
a(\g z)-a(z)  = \int_z^{\g z} \pi^*\alpha = \zeta(\g) - \tfrac{\tau V}{2\pi}S(\g).
$$
We now define 
$$
\triangle_{2\tau}^*(z) \coloneqq  e^{2\pi i\, a(z)}\triangle_{2\tau}(z).
$$
 Then $\triangle_{2\tau}^*$ is a real analytic automorphic form of weight $2\tau$ --- transforming with trivial multiplier system $\chi_{2\tau}^*(\g)=e^{2\pi i\zeta(\g)}$ --- that is nowhere-vanishing on $\h$ and verifies that $\triangle_{2\tau}^*(z)\to0$ as $y\to\infty$. We therefore have a well defined winding number
 ${\rm ind}(\widetilde\triangle^*_{2\tau}(\cdot))$ for homotopy classes of closed curves on the orbifold $M=\G\bk\h$.  We conclude this section with the following direct computation
\begin{align*}
{\rm ind}(\widetilde \triangle^*_{2\tau}(C_\g)) &= \frac{1}{2\pi i} \op{ \log \triangle_{2\tau}(\g z)-\log\triangle_{2\tau}(z) +2\pi ih(\g) -2\tau \log j(\g,z)}\\
& = \frac{1}{2\pi i}\op{ \tau V f(\g z) -\tau V f(z) +2\pi i \zeta(\g) - i\tau V S(\g)-2\tau \log j(\g,z)}\,  = \zeta(\g).
\end{align*}
The definition of $\zeta(\g)$ given above allows to study these winding numbers using combinatorial group theory. We will not pursue this direction here.

\section{Rademacher Symbols for Fuchsian Groups}\label{Rad}

We build on the notation introduced in the previous section. In the literature, Dedekind symbols for Fuchsian groups have been (somewhat confusingly) introduced as one of the following two functions:
$$
S:\G\to\R,\quad S(\g) = \frac{1}{i}\op{ f(\g z) - f(z) -2V^{-1} \log j(\g,z)},
 $$
 where the branch of logarithm is chosen such that $\arg(z)\in(-\pi,\pi]$, or
 $$
\Phi:\G\to\R,\quad \Phi(\g)=\frac{1}{i}\op{ f(\g z)- f(z) - V^{-1}\, \sign(c)^2 \log\op{-j(\g,z)^2}},
$$
where the negative sign is introduced to correct the fact that $\arg (cz+d)^2\in(0,2\pi)$ for $c\neq0$. The choice of branch has the effect that $S(-\g)\neq S(\g)$ if $-I\in\G$, while $\Phi(-\g)=\Phi(\g)$. We record for reference that
 for each $\g_1=\bsm *&*\\ c_1 &d_1\esm$, $\g_2 =\bsm *&*\\ c_2 & d_2\esm$ in $\G$, we have
\begin{align}\label{qm}
\Phi(\g_1\g_2)-\Phi(\g_1) - \Phi(\g_2) = - \pi V^{-1}\sign(c_1 c_2 c_3),
\end{align}
where $\g_1 \g_2 = \bsm *&*\\ c_3 &d_3\esm.$ The reader can further verify that 
\begin{align*}
\Phi(\pm I)& =S(I)=0,\quad S(-I)=-2\pi V^{-1},\\
\Phi(\g^{-1}) &=-\Phi(\g), \quad S(\g^{-1})=-S(\g)
\end{align*}
for each $\g\in\G$. We note that it is the function $\Phi$ that generalizes the Dedekind symbol of Dedekind and Rademacher; compare to Section \ref{subsec:DedRad}. By analogy with the classical Rademacher symbol, the {\bf Rademacher symbol} associated to $\Phi$ is then defined by 
\begin{align}\label{def-psi}
\Psi:\Gamma\to\R,\quad 
\Psi(\gamma) =  \Phi(\gamma) - \pi V^{-1}\sign( c(a+d)),
\end{align}
where $\g=\bsm a& b\\ c& d\esm.$

\begin{lemma}\label{lm-Psi}
By comparison with the properties of the Dedekind symbols $S$ and $\Phi$, we  find that for all $\g\in\G$,
\begin{enumerate}
\item[(i)] $\Psi(-\g)=\Psi(\g)$; 
\item[(ii)] $\Psi(\g^{-1})=-\Psi(\g)$; 
\item[(iii)] $\Psi(\g)=\Phi(\g)=S(\g)$ whenever $\g\in\G_\infty=\{T^n: n\in\Z\}(\times\Z/2\Z)$; 
\item[(iv)] and
$$
\Psi(\g) = \begin{dcases} S(\g) & a+d>0,\\ S(\g) + \pi V^{-1}\sign(c) & a+d=0,\\ S(\g) +2\pi V^{-1}\sign(c) & a+d<0,
\end{dcases}
$$
whenever $\g=\bsm a&b\\ c&d\esm\not\in\G_\infty$.
\end{enumerate}
\end{lemma}

%\begin{proof}
%Statements (i) and (ii) are immediate. For (iii) and (iv), we recall that $\g=\bsm a& b \\ c& d\esm\in\G_\infty$ if and only if $c=0$. This can be seen by looking at the fixed points of $\g$ in $\overline\h=\h\cup\R\cup\{\infty\}$. Indeed, the quadratic equation $\g z=z$ for $\g=\bsm a&b\\ c&d\esm$ has solutions 
%$$
%z_\pm = \frac{a-d\pm \sqrt{(a+d)^2 -4}}{2c}
%$$
%and so $c=0$ corresponds uniquely to the situation in which $z$ is the point at infinity and $\g\in\G_\infty$. Finally, comparing definitions, we find that for $\bsm *&*\\ c&d\esm\in\G$ with $c\neq0$ we have
%\begin{align*}
%\Psi(\g) - S(\g) &= \Phi(\g)-S(\g) - \pi V^{-1}\sign(c(a+d))\\
%& =  -V^{-1} \arg(-(cz+d)^2) +2V^{-1}\arg(cz+d) -\pi V^{-1} \sign(c(a+d))\\
%&= \pi V^{-1} \sign(c)\op{1-\sign(a+d)}.
%\end{align*}
%Statement (iv) follows.
%\end{proof}
 
 One may also write down an algebraic relation of the form  of (\ref{future}) or (\ref{qm}) for $\Psi$ but we omit to do so --- the stated formula for $\Phi$ is usually more convenient to manipulate. Next we prove that $\Psi$ is well defined on the set $\Pi$ of prime geodesics on $M=\G\bk\h$. Since there is bijection between closed geodesics on $M$ and $\G$-conjugacy classes of elements in $\G$, it suffices to show that $\Psi$ is conjugacy class invariant and homogeneous. This follows from the following characterization of the Rademacher symbol on hyperbolic elements as the homogenization of the Dedekind symbol.

\begin{theorem}\label{prop:Psiclassinvariant}
For each $\gamma\in \Gamma$ such that $|\tr(\g)|\geq2$, we have
$$
\Psi(\gamma) = \lim_{n\to\infty} \frac{\Phi(\gamma^n)}{n}. 
$$
\end{theorem}

\begin{proof}
By recursion on \eqref{qm}, we have that for each $n\geq1$,
\begin{align}\label{power}
\Phi(\gamma^n) = n\cdot\Phi(\gamma) -\pi V^{-1}\sum_{k=1}^{n-1} \sign(c_{\gamma} c_{\gamma^k} c_{\gamma^{k+1}}),
\end{align}
where $\g^k = \bsm *&*\\ c_{\g^k}&*\esm$. We will compute $\Phi(\g^n)$ using \eqref{power}. 

Suppose first that $|\mathrm{tr}(\gamma)|> 2$. Then $\g$ has two distinct real eigenvalues $\lambda_1\neq\lambda_2=\tfrac{1}{\lambda_1}$. Up to replacing $\g$ by $\g^{-1}$, we may assume that $\lambda_1>1$.
Then $\gamma^n= g^{-1}\bsm \lambda_1^n & 0 \\ 0 & \lambda_2^{n}\esm g$ for some $g\in G$. Upon comparing matrix coefficients, we find that $$
c_{\gamma^n} = c_\gamma\cdot \frac{\lambda_1^n-\lambda_2^{n}}{\lambda_1 - \lambda_2}.
$$
This implies 
\begin{align*}
\sign(c_\gamma c_{\gamma^n} c_{\gamma^{n+1}}) &= \sign(c_\gamma) \sign((\lambda_1^n-\lambda_2^{n})(\lambda_1^{n+1}-\lambda_2^{n+1}))\\
& = \begin{dcases} \sign(c_\gamma) &\text{ if } \lambda_1>1,\\ -\sign(c_\gamma) & \text{ if } \lambda_1<-1.
\end{dcases}
\end{align*}
The latter conditions being equivalent to $\tr(\gamma)>0$ and $\tr(\gamma)<0$ respectively, we conclude that 
\begin{align}\label{sol}
\lim_{n\to\infty} \frac{\Phi(\gamma^n)}{n} = \Phi(\gamma) -\pi V^{-1}\sign(c(a+d)) = \Psi(\gamma).
\end{align}

Suppose now that $|\tr(\gamma)|=2$. Then $\gamma$ is conjugate (in $G$) to a matrix of the form $\pm\bsm 1&h\\0&1\esm$ for some nonzero $h\neq0$. Hence $\g^n$ is conjugate to $\bsm (\pm 1)^n & (\pm1)^{n-1} nh\\ 0 & (\pm 1)^n\esm$ and comparing matrix coefficients, we now find that $c_{\gamma^n} = c_\gamma (\pm1)^{n-1} n.$ Again, we find $
\sign(c_\gamma c_{\gamma^n} c_{\gamma^{n+1}}) = \sign(c_\gamma \tr(\gamma))$, and consequently, \eqref{sol} holds as well.
\end{proof}

\begin{corollary}\label{conj}
For any $\g\in\G$ with $|\tr(\g)|\geq2$, we have
\begin{itemize}
\item[(i)] $\Psi(\g^n)=n\Psi(\g)$ for each $n\in\Z$,
\item[(ii)] $\Psi(\g_1)=\Psi(\g_2)$ whenever $\g_1,\g_2$ are in the same $\G$-conjugacy class.
\end{itemize}
\end{corollary}

\begin{proof}
Since $\Psi(\g^{-1})=-\Psi(\g)$ and $\Psi(I)=0$, it suffices to proves (i) for $n\in\N$. This follows from
$$
\Psi(\gamma^n) = n\cdot \lim_{m\to\infty}\frac{\Phi(\gamma^{mn})}{mn} = n\cdot \Psi(\gamma).
$$

To prove (ii), we use that by \eqref{qm}, we have the uniform upper bound
\begin{align*}
|\Psi(\tau^{-1}\gamma\tau)-\Psi(\gamma)|\ \leq 4\pi V^{-1},
\end{align*}
for all $\g,\tau\in\G$. Therefore, for any elements $\g_1,\g_2\in\G$ in the same $\G$-conjugacy class, we have
\begin{align*} 
|\Psi(\g_1)-\Psi(\gamma_2)| &=\lim_{n\to\infty} \frac{|\Psi(\gamma_1^n)-\Psi(\gamma_2^n)|}{n} \leq\lim_{n\to\infty} \frac{4\pi V^{-1}}{n}= 0.
\end{align*}
\end{proof}

Given that $\Psi$ is constant on conjugacy classes of hyperbolic elements, the values $\Psi(C)$, with $C\in\Pi$ a prime geodesic, are well defined. Comparing with the construction introduced in the previous section, we can conclude that the Rademacher symbols realize the following winding numbers.

\begin{theorem}[cf.~Theorem \ref{Thm1b}]\label{thm1b}
Let $M$ be a cusped hyperbolic orbifold with finite area $V$. Fix a generating set $\mathcal{S}=\{\g_i\}_{i\in I}$ for $\G$. There exists a nonholomorphic automorphic form $\triangle^*_{2\tau}$ of weight $2\tau$ for $\G$, that is nowhere-vanishing, for which $\triangle^*_{2\tau}(z)\to0$ as $y\to\infty$. Moreover there exists a function $\zeta:\G\to\Z$ such that 
$$
\zeta(\g) = {\rm ind}(\widetilde\triangle_{2\tau}^*(C_\g))
$$
for each hyperbolic element $\g\in\G$ of positive trace, and
$$
\zeta(\g) =\zeta(\g_{i_1}\cdots \g_{i_k}) = \frac{\tau V}{2\pi}\op{\Psi(\g) -\sum_{\substack{ |\tr(\g_{i_j})|\geq2}} \Psi(\g_{i_j})}
$$
for each $\g\in\G$.
\end{theorem}

\begin{theorem}[Theorem \ref{Thm1c}]
Let $M$ be a cusped hyperbolic orbifold with finite area $V$. 
If the Rademacher symbol $\Psi$ is rational-valued, then there exists $k\in2\N$ such that 
$$
{\rm ind}(\widetilde\triangle_k(C))=\frac{kV}{4\pi}\Psi(C)
$$
for each $C\in\Pi$.
\end{theorem}

\section{Estimates on Weyl's Law}\label{sec:Weyl}

In this section we gather estimates from the second author's PhD thesis that will be needed to prove the twisted prime geodesic theorem for multiplier systems. On a first read, the reader is invited to skip directly to Section \ref{sec:PGT} and come back to this section for reference on the relevant estimates.

Evaluating the trace formula reviewed in Section \ref{subsec:STF} for the test function $h(t)=e^{-t^2/T^2}$ yields Weyl's law on the asymptotic distribution of eigenvalues, i.e.,
\begin{align}\label{Weyl}
\sum_{|t_n|\leq T} 1 - \frac{1}{4\pi}\int_{-T}^T \frac{\varphi'}{\varphi}\left(\frac12+it,\chi,r\right)dt\ \sim \frac{V}{4\pi}\ T^2
\end{align}
as $T\to\infty$; see \cite[p.~414]{HejhalVol2}. In general, we have no means of separately estimating the two terms on the left hand-side. (If this were the case, we would obtain from the first term an asymptotic formula for the distribution of the eigenvalues.) We obtain upper bounds on both terms of the LHS of Equation (\ref{Weyl}) as well as their explicit dependence on $\chi$ and $r$. For this purpose, set 
\begin{align}\label{L}
L(\chi,r)\ =\ 1+\sum_{\alpha_\frak{a}(\chi,r)\neq0} \log\op{\alpha_\fa(\chi,r)^{-1}}
\end{align}
where the sum is over all regular cusps $\frak{a}$ of $\Gamma$. The following first result follows by a careful evaluation of the trace formula.

\begin{theorem}\label{wts}
Let $|r|\leq1$ and $T\geq2$. Then
\begin{align*}
\sum_{n\geq0} e^{-t_n^2/T^2} -\frac{1}{4\pi}\int_{-\infty}^\infty e^{-t^2/T^2}\frac{\varphi'}{\varphi}\op{\frac12+it,\chi,r}dt\  \ll\ T^2 + T(L(\chi,r)+1).
\end{align*}
\end{theorem}

\begin{proof}
We evaluate the terms in the trace formula against the test-function $h(t)=e^{-t^2/T^2}$ and its Fourier transform
$$
g(t)\ =\ \frac{T}{2\sqrt\pi} e^{-(tT)^2/4}\ \ll\ T e^{-t^2}.
$$ 
To evaluate \eqref{II}, we note that $|\sinh(t)/\cosh(t)|\to1$ as $|t|\to\infty$, and that $\cosh(t)-1$ is positive except in $t=0$, where it has a double zero. Since $\sinh(t)$ has a zero at $t=0$, the integrand
$$
\frac{t\sinh(2\pi t)}{\cosh(2\pi t)+\cos(\pi r)}
$$
can be continuously extended to $(t,r)=(0,\pm1)$ and is uniformly bounded for $|t|\leq1$. Hence
\begin{align*}
\eqref{II} &\ll\  \int_{-\infty}^\infty (|t|+1)e^{-t^2/T^2}  dt\ \ll\ T^2.
\end{align*}
Since we only consider $|r|\leq1$, \eqref{I} is 0. For \eqref{E}, there are only finitely many terms, and
\begin{align*}
\eqref{E}\ &\ll\  T \int_{-\infty}^\infty e^{-t^2+(r-1)t/2} (e^t +1) dt\ \ll\ T.
\end{align*}
The same bound also applies to the terms (\ref{B}) and (\ref{C}). We conclude that \begin{align*}
|\eqref{A}+\eqref{B}+\eqref{C}+\eqref{D}|\ & \ll\ T( L(\chi,r) +1) +\int_{-\infty}^\infty e^{-t^2/T^2} \left|\frac{\Gamma'}{\Gamma}(1+it)\right| dt.
\end{align*}
Further, we rely on the standard approximation \cite[(B.11)]{Iwaniec}
\begin{align*}
\frac{\Gamma'}{\Gamma}(1+it)\ =\ \log(s) -\frac{1}{2s} +O\left(\frac{1}{|s|^2}\right)
\end{align*}
to bound
\begin{align*}
\int_{-\infty}^\infty e^{-t^2/T^2} \left|\frac{\Gamma'}{\Gamma}(1+it)\right| dt\ &\ll\ \int_{-\infty}^\infty e^{-t^2/T^2}(|t|+1)dt\ \ll\ T^2.
\end{align*}
For the hyperbolic contribution, we have

\begin{align*}
\eqref{H}\ &\ll  T\sum_{\substack{\{\gamma\}_{\rm pr}\\ \tr(\gamma)>2}} \frac{\ell_\g}{\sinh(\ell_\g/2)}\sum_{k\geq1} e^{-k\ell_\g^2} \ll T\sum_{\substack{\{\gamma\}_{\rm pr}\\ \tr(\gamma)>2}} \frac{\ell_\g e^{-\ell_\g}}{\sinh(\ell_\g/2)}\, =\, O(T),
\end{align*}
where the last equality follows from the prime geodesic theorem. Finally, using that the entries in the scattering matrix $\Phi(\tfrac12,\chi,r)$ are bounded and combining these estimates, we conclude that Theorem \ref{wts} holds.
\end{proof}

We show next how to obtain separate upper bounds for the two terms on the LHS of (\ref{Weyl}) using the Dirichlet series representation (\ref{Dirich}) and a positivity argument.

\begin{theorem}\label{thmFvE}
If $|r|\leq1$, we have
\begin{align}\label{disc}
|\{ n\in\N: |t_n|\leq T\}|\
 \ll T^2 + T(L(\chi,r)+1),
\end{align}
and
\begin{align}\label{cont}
\int_{-T}^T \left|\frac{\varphi'}{\varphi}\left(\frac12+it,\chi,r\right)\right|dt\ \ll\ T^2 + T(L(\chi,r)+1).
\end{align}
\end{theorem}

\begin{proof}
We consider the function
\begin{align}\label{Selberg}
\varphi^*(s,\chi,r) = b_1^{2s-1}\prod_{j=1}^{k}\frac{s_j(\chi,r)-s}{s_j(\chi,r)+s-1}\varphi(s,\chi,r),
\end{align}
where $b_1>0$ is the smallest positive denominator in the Dirichlet series given  by \eqref{Dirich} and is bounded below. The local factors are repeated with multiplicity for each pole $s_j(\chi,r)\in(1/2,1]$ --- see Lemma \ref{lm:mult} ---, and $\varphi^*(s,\chi,r)$ is regular and uniformly bounded for $\re(s)\geq1/2$. Selberg \cite[pp.~655-656]{Selberg1989} shows that 
$$
-\frac{\varphi^{*}~'}{\varphi^*}\left(\frac12+it,\chi,r\right)  >0
$$
for all $t\in\R$. 
By taking the logarithmic derivative of \eqref{Selberg}, we have
\begin{align}\label{relation}
-\frac{\varphi^{*}~'}{\varphi^*}(s,\chi,r) = -\frac{\varphi'}{\varphi}(s,\chi,r) - 2\log b_1 + \sum_{j=1}^k \frac{2s_j(\chi,r)-1}{(s_j(\chi,r)-s)(s_j(\chi,r)+s-1)}.
\end{align}
Since $b_1$ is uniformly bounded below, the term $-\log b_1$ is uniformly bounded above. For $|t|\leq T$, $e^{-t^2/T^2}\geq e^{-1}$, and using that $-\tfrac{\varphi^{*}~'}{\varphi^*}>0$, we have 
\begin{align*}
\sum_{|t_n|\leq T} 1 -\frac{1}{4\pi}\int_{-T}^T\frac{\varphi^{*}~'}{\varphi^*}\left(\frac12+it,\chi,r\right) dt& \ll \sum_{|t_n(r)|\leq T} e^{-\tfrac{t_n^2}{T^2}}\\
&\quad -\frac{1}{4\pi}\int_{-T}^T e^{-\tfrac{t^2}{T^2}}\frac{\varphi^{*}~'}{\varphi^*}\left(\frac12+it,\chi,r\right) dt.
\end{align*}
Using the relation \eqref{relation}, the last line implies the estimate
\begin{align*}
\sum_{|t_n|\leq T} 1 - &\frac{1}{4\pi}\int_{-T}^T\frac{\varphi'}{\varphi}\left(\frac12+it,\chi,r\right) dt \\
&\qquad  \ll\
\sum_{|t_n|\geq0} e^{-t_n(r)^2/T^2} -\frac{1}{4\pi}\int_{-\infty}^\infty e^{-t^2/T^2}\frac{\varphi'}{\varphi}\left(\frac12+it,\chi,r\right) dt + E(T)
\end{align*}
for 
$$
E(T) = \sum_{j=1}^k \frac{2s_j(\chi,r)-1}{4\pi} \int_{-T}^T \frac{ e^{-t^2/T^2}-1}{(\sigma_j(r)-\tfrac12)^2+t^2}dt -\frac{\log b_1}{2\pi}\left(2T-\int_{-T}^T e^{-t^2/T^2} dt\right).
$$
We have $E(T)=O(T)$. Theorem \ref{thmFvE} then follows from Theorem \ref{wts} and evaluating the terms individually using the positivity of $-\tfrac{\varphi^{*}~'}{\varphi^*}$.
\end{proof}

\section{Prime Geodesic Theorem for Multiplier Systems}\label{sec:PGT}

We prove Theorem \ref{Thm2}.

\begin{theorem}[Theorem \ref{Thm2}]
Let $\G$ be a cofinite Fuchsian group and let $\chi$ be a multiplier system of weight $r$ on $\G$. Then as $T\to\infty$, we have
\begin{align}
\sum_{C_\g\in \Pi(T)} \chi(\g) =  \Li\left(e^{s_0(\chi,r)T}\right)+\dots+\Li\left(e^{s_k(\chi,r)T}\right) + O\left(e^{3T/4}L(\chi,r)\right),
\end{align}
where $s_0(\chi,r)\geq s_1(\chi,r) \geq \dots\geq s_{k}(\chi,r)>\tfrac12$ are the spectral eigenparameters determined by the small eigenvalues $\lambda_j(\chi,r)=s_j(\chi,r)(1-s_j(\chi,r))$ of the weight $r$ Laplacian $\Delta_r$, and $L(\chi,r)$ is defined by (\ref{L}).
\end{theorem}

\subsection{Choice of test-function.} We choose $g$ to be the mollified characteristic function on the symmetric interval $[-T,T]$. More precisely, we choose a smooth even function $k\in C_c^\infty(\R)$ of compact support such that $\int_\R k(t)dt=1$ and $k(t)\geq0$. Then fix $\eps>0$, let $k_\eps(t)=\tfrac{1}{\eps}k(\tfrac{t}{\eps})$, and set up the convolution product
$$
g(t)=(\mathbf{1}_{[-T,T]}\ast k_\eps)(t)=\int_{-T}^T k_\eps(t-u)du.
$$
Observe that $\eps$ controls how closely $g$ approximates $\mathbf{1}_{[-T,T]}$; indeed, we have
$$
g(t) =\begin{dcases} 1 & \text{ if } |t|\leq T-\eps,\\ 0 & \text{ if } |t|\geq T+\eps,\end{dcases}
$$
and $g$ is smoothly decaying on $(T-\eps,T+\eps)$. The parameter $\eps$ will be chosen later on as a function of $T$. The inverse Fourier transform $h$ is given by
$$
h(t)\ =\ \frac{2\sin(tT)}{t} \widehat{k}(\eps t).
$$
We again underline that both $g$ and $h$ depend on the parameters $T, \eps$. In particular, $g(0)=1$ if $\eps\leq T$ and $h(0)=2T\widehat{k}(0)=2T$. By integration by parts, we have
$$
|\widehat{k}(\eps t)|\leq \frac{1}{(\eps t)^j}
$$
for any $j\geq1$. As $|t|\to0$, we have 
$$
\widehat{k}(\eps t) = 1 + O(\eps|t|).
$$

\subsection{Basic estimates.} 
We now collect estimates for the trace formula. We have
\begin{align*}
\eqref{II} & \ll\ \int_{-\infty}^\infty |h(t)|(|t|+1)dt \ll T\int_{-\infty}^\infty (|t|+1)\frac{|\sin(tT)|}{|tT|}|\widehat{k}(\eps t)|dt\\
& \ll\ T\int_{-1}^1 \left|\frac{\sin(tT)}{tT}\widehat{k}(\eps t)\right|dt + 2\int_1^\infty \left|\sin(tT)\widehat{k}(\eps t)\right| dt\\
&\ll T + \eps^{-1} \int_1^\infty |\widehat{k}(t)|dt\ \ll\ T + \eps^{-1},\\
\eqref{E}\ & \ll\ \int_{-\infty}^\infty |g(t)| dt\ \ll\ T,
\end{align*}
where we used that the sum is finite and that for fixed $\theta$, the term
$$
e^{(r-1)t/2}\frac{e^t-e^{2i\theta}}{\cosh(t)-\cos(2\theta)}
$$
is uniformly bounded. Similarly, (\ref{B}) and (\ref{C}) are $O(T)$. For (\ref{D}), we have
\begin{align*}
\int_{-\infty}^\infty |h(t)|\left|\frac{\Gamma'}{\Gamma}(1+it)\right|dt &\ll
\int_{-\infty}^\infty |h(t)\log(1+it)|dt + \int_{-1}^1 |h(t)|dt\\
&\ll T + \int_{-\infty}^\infty |\widehat{k}(\eps t)| dt\\
& \ll T + \eps^{-1}\int_{-\infty}^\infty (1+t^2)^{-1} dt \ll T+\eps^{-1}.
\end{align*}
 For the hyperbolic contribution, we use the estimate 
$$
\eqref{H} = \sum_{C_\g\in\Pi(T)} \frac{\chi(\gamma)\ell_\gamma}{\sinh(\ell_\gamma/2)} + O(\eps e^{T/2} +T^2).
$$
Finally, using Theorem \ref{thmFvE}, we have
\begin{align*}
\eqref{S}\ &\ll\ m(\chi,r) T + \int_{-\infty}^\infty |h(t)|\left|\frac{\varphi'}{\varphi}\left(\frac12+it,\chi,r\right)\right|dt\\
& \ll m(\chi,r)T + \int_{-1}^1 (1+O(\eps|t|))\left|\frac{\varphi'}{\varphi}\left(\frac12+it,\chi,r\right)\right|dt\\
&\qquad  + \int_1^{\eps^{-1}} \frac{1}{t}\left|\frac{\varphi'}{\varphi}\left(\frac12+it,\chi,r\right)\right|dt +\eps^{-2} \int_{\eps^{-1}}^\infty \frac{1}{t^3}\left|\frac{\varphi'}{\varphi}\left(\frac12+it,\chi,r\right)\right|dt\\ 
&\ll (m(\chi,r)+L(\chi,r)+1)T + \eps^{-1}L(\chi,r).
\end{align*}
Collecting estimates, the trace formula yields
$$
\sum_{n\geq0} h(t_n) = \sum_{C_\g\in\Pi(T)} \frac{\chi(\g)\ell_\g}{\sinh(\ell_\g/2)} + O(\eps e^{T/2} + T^2 + \eps^{-1}L(\chi,r)).
$$

\subsection{Small eigenvalues.} 

For each eigenvalue $\lambda_j=\lambda_{j}(\chi,r)<\tfrac14$, the parametrization $\lambda_{j}=(\tfrac12+it_j)(\tfrac12-it_j)=\tfrac14+t_{j}^2$ implies that we may choose $t_{j}\in -i(0,\tfrac12]$. Then
$$
h(t_{j})\ =\ \frac{e^{i t_{j}T}-e^{-i t_{j} T}}{i t_{j}(r)}\widehat{k}(\eps t_{j})\ =\ \frac{e^{i t_{j} T}}{i t_{j}} + O(\eps e^{|t_{j}|T})
$$
and hence
$$
\sum_{\lambda_{j}(\chi,r)<1/4} h(t_{j})\ =\  \sum_j \frac{e^{i t_{j}T}}{i t_{j}} + O(\eps e^{T/2}),
$$
where the implied constant is independent of $\chi$ and $r$. 

\subsection{Cuspidal eigenvalues.} Consider now $\lambda_{j}(\chi,r)\geq1/4$, and hence $t_{j}\in\R$. Then
$$
\sum_{0\leq|t_{j}|<1} h(t_{j}) =  2T\sum_{0\leq|t_{j}|<1} \frac{\sin(t_{j}T)}{t_{j} T} \widehat{k}(\eps t_{j}) = O(T),
$$
while by summation by parts and Theorem \ref{thmFvE}, we have
$$
\sum_{1\leq |t_{j}|<\eps^{-1}} h(t_{j})\ \ll\ \sum_{1\leq|t_{j}|<\eps^{-1}} \frac{1}{|t_{j}|}\ \ll\ \eps^{-1} L(\chi,r),
$$
and using additionally the fast decay of Fourier coefficients, we have
$$
\sum_{1\leq\eps|t_{j}|} h(t_{j})\ \ll\ \eps \sum_{1\leq\eps|t_{j}|} \frac{1}{(\eps|t_{j}|)^3}\ \ll\ \eps^{-1}L(\chi,r).
$$

\subsection{Bootstrapping.}
Collecting estimates, we conclude that 
$$
\sum_{C_\g\in\Pi(T)} \frac{\chi(\gamma)\ell_\gamma}{\sinh(\ell_\gamma/2)}\ =\  \sum_j \frac{e^{i t_{j} T}}{i t_{j}} +
 O(\eps e^{T/2}+T^2+\eps^{-1}L(\chi,r)),
$$
where the error term can be optimized to $O(e^{T/4}(L(\chi,r)+1))$. Applying summation by parts  we find 
\begin{align}
\sum_{C_\g\in\Pi(T)} \chi(\gamma)\ell_\gamma\ &=\ \sum_{j} \frac{e^{(i t_j+1/2)T}}{i t_j+1/2} + O(e^{3T/4}L(\chi,r))\nonumber\\
&=\ \frac{e^{s_0(\chi,r)T}}{s_0(\chi,r)} + \frac{e^{s_1(\chi,r)T}}{s_1(\chi,r)} + \dots + \frac{e^{s_k(\chi,r)T}}{s_k(\chi,r)} + O\left(e^{3T/4} L(\chi,r)\right)\label{stepstone}
\end{align}
and 
\begin{align*}
\sum_{C_\g\in\Pi(T)} \chi(\g) &= \Li\op{e^{s_0(\chi,r)T}} +\dots + \Li\op{e^{s_k(\chi,r)T}} + O\op{e^{3T/4}L(\chi,r)},
\end{align*}
which concludes the proof of Theorem \ref{Thm2}.

\section{Small Eigenvalues and Continuous Deformations}\label{sec:stab}

We consider a continuous perturbation of (\ref{stepstone}) in the weight $r$. For this we need to understand how the small eigenvalues behave under perturbation. The main result of this section is 
\begin{lemma}\label{lm2}
Let $\G$ be a cofinite Fuchsian group and let $\chi$ be a multiplier system of weight $r$ for $\G$ such that $\chi(\g)$ is continuous as a function of $r$, for any fixed $\g\in\G$. Fix $T<1/4$ such that it is not an eigenvalue for $\Delta_0=\Delta$. Then the count 
$\#\{\lambda_{n}(\chi,r)\leq T\}$ is continuous at $r=0$.
\end{lemma}

Our proof follows the strategy of \cite[Lemma 3.3]{Risager2011} (which proves the same result for automorphic forms transforming with a unitary character) together with an estimate from \cite{JorgensonLundelius1997} for the hyperbolic heat trace. We will need the following facts pertaining to the Laplace transform. For a sufficiently nice function $f:\R_{>0}\to\C$, its Laplace transform is given by
$$
\mathscr{L}(f)(z) = \int_0^\infty e^{-zt} f(t) dt
$$
with inverse transform
$$
f(u) = \frac{1}{2\pi i}\int_{(a)} e^{zu} \mathscr{L}(f)(z)dz.
$$
Let $\rho>0$ and define
$$
f_\rho(t) = \int_0^t \frac{(t-u)^{\rho-1}}{\Gamma(\rho)}f(u) du.
$$
Then following \cite[Thm. 8.1, 8.2, p.~73]{Widder1941},  we have
\begin{align}\label{Laplace}
\frac{1}{2\pi i}\int_{(a)} e^{zu} \frac{\mathscr{L}(f)(z)}{z^\rho}dz = \begin{dcases}
f_\rho(u) & \text{ if } u\geq0,\\
0 & \text{ if } u<0,
\end{dcases}
\end{align}
whenever $a>0$ is large enough for $\mathscr{L}(f)(z)$ to converge absolutely in the half-plane $\re(z)\geq a>a_0$.

\begin{proof}[Proof of Lemma \ref{lm2}]
Let $f(t)=t^{w-1}$ with $w\geq1$. Its Laplace transform is 
$$
\mathscr{L}(f)(z) = \frac{\Gamma(w)}{z^w}.
$$
We consider the trace formula for the test function $h_z(t)=e^{-zt^2}$ for a fixed $z\in\C$ with fixed positive real part $a=\re(z)>0$ and we integrate the trace formula against 
\begin{align}\label{integrate}
\frac{\mathscr{L}(f)(z) e^{z(T-1/4)}}{2\pi i z} = \frac{\Gamma(w) e^{z(T-1/4)}}{2\pi i z^{w+1}}
\end{align}
along $(a)=a+i\R$. On the left hand-side of the trace formula, this yields with (\ref{Laplace})
\begin{align*}
\frac{1}{2\pi i} \sum_{n\geq0} \int_{(a)} e^{z(T-(\tfrac14+t_n^2))} \frac{\mathscr{L}(f)(z)}{z}dz\ & =\ 
\sum_{\lambda_n(\chi,r)\leq T} f_1(T-\lambda_n(r))\\
& =\ \sum_{\lambda_n(\chi,r)\leq T} \int_0^{T-\lambda_n(\chi,r)} f(u)du\\
& =\ \frac{1}{w} \sum_{\lambda_n(\chi,r)\leq T} (T-\lambda_n(\chi,r))^w.
\end{align*}
Recall that $a=\re(z)>0$. The Fourier transform of $h_z(t)$ is given by
$$
g_z(t) = \sqrt{\frac{\pi}{z}} e^{-t^2/(4z)}.
$$
The integration against (\ref{integrate}) annihilates most terms of the trace formula. Indeed, for all $t\in\R$ and $T<\tfrac14$, we have
$$
\frac{1}{2\pi i} \int_{(a)} h_z(t) \frac{e^{z(T-1/4)}}{z^2}dz = \frac{1}{2\pi i}\int_{(a)} \frac{e^{z(T-\tfrac14-t^2)}}{z^2}dz =0,
$$
and
$$
\frac{1}{2\pi i} \int_{(a)} g_z(0) \frac{e^{z(T-1/4)}}{z^2}dz = \frac{1}{2\sqrt{\pi} i}\int_{(a)} \frac{e^{z(T-\tfrac14)}}{z^{5/2}} dz =0.
$$
By (6.43)--(6.44) in \cite[p.401]{HejhalVol2}, the integrals given by (\ref{E}) and (\ref{B}) can be expressed in terms of $h_z(t)$ rather than $g_z(t)$. So these terms vanish as well after integration. We are left with
\begin{align*}
N_w(T,r) := \frac{1}{w}\sum_{\lambda_n(\chi,r)\leq T} (T-\lambda_n(r))^w\ =\ \frac{1}{2\pi i}\int_{(a)} F_r(z)\frac{\Gamma(w) e^{z(T-1/4)}}{z^{w+1}}dz,
\end{align*}
where 
$$
F_r(z) = \sum_{\{\g\}_{\rm pr}} \sum_{k\geq1} \frac{\chi(\g)\ell_\g g_z(k\ell)}{e^{k\ell/2}-e^{-k\ell/2}} + m(\chi,r)\int_0^\infty \frac{g_z(u)(1-\cosh(\tfrac{ru}{2}))}{e^{u/2}-e^{-u/2}}du.
$$

We first claim that $N_w(T,r)$ is continuous at $r=0$ by dominated convergence for $w>\tfrac32$. First, we see that, as $r\to0$,
$$
\frac{g_z(u)(1-\cosh(\tfrac{r}{2}u))}{e^{u/2}-e^{-u/2}} \frac{e^{z(T-1/4)}}{z^{w+1}}\ \to\ 0
$$
uniformly for $z\in (a)$, and $t\in\R_{>0}$. By dominated convergence, the second term of $F_r(z)$ is then continuous as a function of $r$ at $r=0$, and in fact 0 at $r=0$. (In a small interval $|r|\leq\delta$, $r\neq0$, the count $m(\chi,r)$ is constant.) To control the hyperbolic contribution, we use the bound 
$$
|H(a+it)| \ll (1+|t|)^{3/2}.
$$
(see (4.2) in \cite{JorgensonLundelius1997}) for the hyperbolic heat trace
$$
H(z) = \frac{e^{-z/4}}{4\sqrt{\pi z}} \sum_{\{\g\}_{\rm pr}} \sum_{k\geq1} \frac{\ell_\g e^{-(k\ell_\g)^2/(4z)}}{\sinh(k\ell_\g/2)}.
$$
Then
$$
\frac{1}{2\pi i} \int_{(a)} H(z) \frac{e^{z(T-1/4)}}{z^{w+1}} dz
$$
is bounded by
$$
\int_1^\infty \frac{(1+t)^{3/2}}{(a^2+t^2)^{\tfrac{1+w}{2}}} dt + O(1) \ll  \int_1^\infty t^{1/2-w} dt.
$$
This proves that $N_w(T,r)$ is continuous at $r=0$ for $w>3/2$. To extend the claim to all $w\geq1$, we first observe that by definition, the function $N_w(T,r)$ is positive and monotone increasing in $T$ and that 
\begin{align}\label{use2}
\frac{d}{dT} N_{w+1}(T,r) = w N_w(T,r),
\end{align}
where we have used that $\G(w+1)=w\G(w)$. Hence by the mean value theorem, we have
$$
N_w(T,r) \leq \frac{N_{w+1}(T+\delta,r)-N_{w+1}(T,r)}{w\delta} \leq N_w(T+\delta,r)
$$
for $T<T+\delta<1/4$. So for any $w\geq1$, using the continuity of $N_{w+1}$ at $r=0$, we have
$$
\limsup_{r\to0} N_w(T,r)\ \leq\ \frac{N_{w+1}(T+\delta,0)-N_{w+1}(T,0)}{w\delta}\ =\ N_w(T,0),
$$
where the last equality follows from (\ref{use2}). Similarly, we have
$$
\liminf_{r\to0} N_w(T,r) \geq \frac{1}{w+1} \frac{N_{w+1}(T,0)-N_{w+1}(T-\delta,0)}{\delta} = N_w(T,0),
$$
and hence $N_w(T,r)$ is continuous at $r=0$ for all $w\geq1$ and in particular for $w=1$. By positivity, we have
$$
\frac{1}{\delta}(N_1(T,r)-N_1(T-\delta,r))\ \leq\ \#\{\lambda_{n}(r)\leq T\}\ \leq\ \frac{1}{\delta}(N_1(T+\delta,r)-N_1(T,r))
$$
for $0<T+\delta<\tfrac14$. It follows that the counting function in the middle is continuous at $r=0$. 
\end{proof}

\begin{corollary}\label{const}
Fix $T<1/4$ that is not an eigenvalue for $\Delta_0=\Delta$. There exists $\delta>0$ such that for $|r|\leq\delta$, the number $\#\{\lambda_n(\chi,r)\leq T\}$ of small eigenvalues counted with multiplicity is constant.
\end{corollary}

Recall that $\lambda_0(\chi,0)=0$ for $L^2(\Gamma\backslash\h,\chi,0)$ if and only if $\chi$ is trivial at $r=0$. 
\begin{corollary}\label{conclu}
If $\chi$ is trivial at $r=0$, then the bottom eigenvalue $\lambda_0(\chi,r)$ has multiplicity 1 in some small interval $|r|\leq \delta$.
\end{corollary}

\section{Statistics for Winding Numbers}\label{sec:stat}

Let 
$$
u_0(z) \coloneqq \begin{dcases} y^{r/2} \triangle_r(z) & \text{if }r\geq0,\\
y^{-r/2}\overline{\triangle_r(z)} & \text{if }r\leq 0,
\end{dcases}
$$
where $\triangle_r$ is defined as in (\ref{def:triangle}). 
The following observation will play an important role in  our results.
\begin{proposition}\label{prop:base ef}
We have $u_0\in L^2(\G\bk\h,\chi_r,r)$ and 
$$
\Delta_r u_0 =  \frac{|r|}{2}\op{1-\frac{|r|}{2}}u_0,
$$
where $\Delta_r$ is the weight $r$ Laplacian.
\end{proposition}
\begin{proof}
The transformation formula for $\triangle_r$ implies that $u_0$ satisfies \eqref{real wt}, while the definition of $f$ guarantees that $\scal{u_0,u_0}<\infty$. A direct computation shows that for $u_0=y^{r/2}\triangle_r$, we have
$$
\Lambda_r u_0 = y^{1+r/2}\op{i \frac{\partial}{\partial x}\triangle_r -  \frac{\partial}{\partial y}\triangle_r} 
$$
while for $u_0=y^{-r/2}\overline\triangle_r$, we have
$$
K_r u_0 = y^{1-r/2}\op{i\frac{\partial}{\partial x}\overline\triangle_r -  \frac{\partial}{\partial y}\overline\triangle_r}.
$$
Since $\triangle_r$ is holomorphic, we can conclude that $\Lambda_r u_0 = K_r u_0=0$.
\end{proof}

 \begin{theorem}[Theorem \ref{Thm1}]
 Let $\G\in\mathcal{G}$ and consider the winding number given by Theorem \ref{Thm1c}. Then $\pi_n(T)$ is the number of prime geodesics in $\Pi(T)$ with winding 
 $$
 {\rm ind}(\widetilde\triangle_k(C))=n
 $$ 
 and there exists $\delta\in(0,\tfrac12]$ such that 
$$
\pi_n(T) = \frac{4}{k T}\int_2^{e^T} \frac{\log t}{(\log t)^2 +\op{\tfrac{4\pi n}{k}}^2}\, dt\, + O\op{\frac{e^{T(1-\delta/2)}}{T}}
$$
as $T\to\infty$. The implied constant does not depend on $n$.
\end{theorem}

\begin{proof}
We specialize the twisted Prime Geodesic Theorem \ref{Thm2} to the multiplier system $\chi_r$ associated to $\triangle_r$ in (\ref{def:triangle}). Since $\G\in\mathcal{G}$, Theorem \ref{Thm1c} guarantees the existence of a weight $k\in2\N$ such that $\chi_k$ precisely encodes the winding number
$$
\frac{kV}{4\pi}\Psi(\g) = {\rm ind}(\widetilde\triangle_k(C_\g))
$$
in the following way. Choosing for each $C_\g$ a representative $\g\in\G$ with positive trace, we have $$\chi_r(\g)=e^{2\pi i\frac{rV}{4\pi}\Psi}=e^{2\pi i \frac{r}{k}{\rm ind}(\widetilde\triangle_k(C_\g))}.
$$ 
By periodicity, we may restrict $r$ to the interval $I=(-k/2,k/2]$. We have
\begin{align*}
\sum_{C\in\Pi(T)} \ell_C \int_{I} e^{ 2\pi i \tfrac{r}{k}\op{{\rm ind}(\widetilde\triangle_k(C))-n}}dr\ &=\ k\sum_{\substack{C\in \Pi(T)\\ {\rm ind}(\widetilde\triangle_k(C))=n}} \ell_C.
\end{align*}
Integrating on the RHS of (\ref{stepstone}) accordingly will require more care.

First consider the case of $|r|<1$.  By Proposition \ref{prop:base ef},  we have $\lambda_0(\chi_r,r)=\tfrac{|r|}{2}(1-\tfrac{|r|}{2})$ or, equivalently, $s_0(\chi_r,r)=1-\tfrac{|r|}{2}$. In particular, $\lambda_0(\chi_0,0)=0$ with multiplicity 1 and we may thus choose $T$ small enough so that Lemma \ref{lm2} asserts the existence of $\delta>0$ such that $N_r(T)=N_0(T)=1$ whenever $|r|<\delta$. That is, $\lambda_1(\chi_r,r)\geq T$ whenever $|r|<\delta$, or equivalently, there is a constant $c$ (depending only on $T$) such that $s_1(\chi_r,r)\leq 1-\tfrac{c}{2}$ whenever $|r|<\delta$. Up to choosing a smaller $\delta$, we may assume that $\delta\leq \min\{c,\tfrac12\}$. Then whenever $|r|<\delta$, the twisted prime geodesic theorem in the form given by (\ref{stepstone}) yields
\begin{align}\label{theLHS}
\sum_{C_\g\in \Pi(T)} \chi_r(\g)\ell_\g = \frac{e^{T(1-\tfrac{|r|}{2})}}{1-\tfrac{|r|}{2}} + O\left( e^{T(1-\tfrac{\delta}{2})} L(\chi_r,r)\right)
\end{align}
whereas the LHS is trivially bounded by the error term alone when $\delta\leq |r|\leq1$, since $s_0(\chi_r,r)\leq1-\tfrac{|r|}{2}\leq1-\tfrac{\delta}{2}$.

\vspace{.2cm}
If $|r|>1$, we fix $\rho\in(-1,1]$ to be adjusted weight uniquely determined by $\rho\equiv r$ (mod 2); that is, $\chi_r = \chi_\rho \chi$ for some nontrivial unitary character $\chi$. Since then $\chi_r=\chi\neq1$ at $\rho=0$, we have $\lambda_1(\chi,0)\geq \lambda_0(\chi,0)>0$ and again there exists $\delta>0$ such that $\lambda_0(\chi_r,r)\geq T$ whenever $|\rho|<\delta$, and we conclude as above that the LHS in (\ref{theLHS}) is $O(e^{T(1-\delta/2)})$ for $|r|> 1$. This proves that there exists $\delta\in(0,1/2]$ such that
\begin{align}\label{thm:refined}
\sum_{C_\g\in\Pi(T)} \chi_r(\g) \ell_\g\ =\ \begin{dcases}
\frac{e^{T(1-\tfrac{|r|}{2})}}{1-\tfrac{|r|}{2}} + O\left(e^{T(1-\tfrac{\delta}{2})}L(\chi_r,r)\right) & \text{if } |r|\leq \delta,\\
O\left(e^{T(1-\tfrac{\delta}{2})}L(\chi_r,r)\right) & \text{if }|r|>\delta.
\end{dcases}
\end{align}

\vspace{.2cm}

 Now integrating the RHS of (\ref{thm:refined}) against $e^{-2\pi i \tfrac{rn}{k}}dr$ over the interval $I$ yields
 \begin{align*}
 \int_{-\delta}^\delta \frac{e^{T(1-|r|/2)}}{1-|r|/2}e^{- 2\pi i \tfrac{rn}{k}}\, dr\
& =\ 2\re \int_0^\delta \frac{e^{T(1-r/2)}}{1-r/2}e^{- 2\pi i \tfrac{rn}{k}}\, dr\\
&=\ 2\re \int_0^\delta \left( \int_2^{e^T} \frac{dy}{y^{r/2}}\right) e^{- 2\pi i \tfrac{rn}{k}}\, dr + O(1)\\
& = 2\re \int_2^{e^T} \int_0^\delta e^{-\tfrac{r}{2}\op{\log y -4\pi i \tfrac{n}{k}}} dr\ dy +O(1)\\
&= 4 \int_2^{e^T} \frac{\log y}{(\log y)^2+(\tfrac{4\pi}{k}n)^2}\, dy + O\left(\int_2^{e^T}  \frac{e^{-\tfrac{\delta}{2} \log y}}{\log y}\, dy\right)\\
& = 4\, \int_2^{e^T} \frac{\log y}{(\log y)^2+(\tfrac{4\pi}{k}n)^2}\, dy+  O\left(\frac{e^{T(1-\delta/2)}}{T}\right).
\end{align*}

It remains to integrate the error term in Theorem \ref{Thm2}. Using \cite[Proposition 5.6]{JorgensonOSullivanSmajlovic2022} we can check that each cusp of $\G$ is singular except the cusp at infinity, for which we have $\alpha_\infty(\chi_r,r)=\tfrac{rV}{4\pi}$.  
Then for every $\eps>0$ we have
$$
\int_{\eps}^1 L(\chi_r,r)\, dr = - \int_\eps^1 \log(r)\, dr + O(1) =  O(1)
$$
as $\eps\to0^+$. We conclude with an application of summation by parts.
\end{proof}

\begin{theorem}[Theorem \ref{Thm:Cauchy}]
Let $\G\in\mathcal{G}$ and consider the winding number given by Theorem \ref{Thm1c}. 
 The limiting distribution of the ratio of winding-to-length for prime geodesics is Cauchy. More precisely, for each $t\in\R$ we have
$$
\lim_{T\to\infty} \frac{\#\{ C\in \Pi(T) : {\rm ind}(\widetilde\triangle_k(C))\leq \tfrac{k\, t}{4\pi}\ell_C\}}{\pi(T)} = \int_{-\infty}^{t} \frac{du}{\pi(1+u^2)}.
$$
\end{theorem}

\begin{proof}
An application of integration by parts yields
\begin{align}\label{IBP}
\pi_n(T) = \frac{4}{kT}\int_2^{e^T} \frac{\log y}{(\log y)^2+(\tfrac{4\pi}{k}n)^2}\, dy =
\frac{4}{k} \frac{Te^T}{T^2+(\tfrac{4\pi}{k}n)^2} +O\op{\frac{e^T}{T^3}}.
\end{align}
We see the Cauchy distribution appear in
\begin{align*}
\sum_{aT\leq n\leq bT} \pi_n(T)\ &=\ \frac{4}{k}\frac{e^T}{T^2}\op{\sum_{aT\leq n\leq bT} \frac{1}{1+\op{\tfrac{4\pi}{k T}n}^2} + O(1)}\\
&= \frac{4}{k}\frac{e^T}{T^2}\op{ \int_{aT}^{bT}  \frac{dx}{1+\op{\tfrac{4\pi}{k T}x}^2}+O(1)}\\
& = \frac{e^T}{T} \int_{4\pi a/k}^{4\pi b/k} \frac{du}{\pi(1+u^2)} + O\op{\frac{e^T}{T^2}},
\end{align*}
where we have used Euler--Maclaurin and a change of variable. For simplicity, we may assume that $b>0$. Let $\eps>0$. Then
\begin{align*}
\frac{1}{\pi(T)} \sum_{\substack{C\in\Pi(T)\\ b\ell_C \leq{\rm ind}(\widetilde\triangle_k(C))\leq bT}} 1\,
& \leq\, \frac{\pi(T(1-\eps))}{\pi(T)} + \frac{1}{\pi(T)}\sum_{bT(1-\eps)\leq n \leq bT} \pi_n(T).
\end{align*}
Since this is true for any $\eps>0$, we conclude that the contribution is 0 in the limit $T\to\infty$. 
\end{proof}

\begin{definition}
A function $f:\Gamma\to\R$ is a {\bf quasimorphism} if there exists $D>0$ such that
$$
\sup_{\gamma_1,\gamma_2\in\Gamma} |f(\gamma_1\gamma_2)-f(\gamma_1)-f(\gamma_2)|\ \leq\ D\ < +\infty.
$$
The constant $D$ is called the defect of the quasimorphism $f$.
\end{definition}
From the properties proved in Section \ref{Rad}, we can see that the Dedekind and Rademacher symbols are examples of quasimorphisms. To support our heuristic that the appearance of the Cauchy law reflects the behavior of long geodesics winding high in the cusp, we observe that conjugacy class invariant quasimorphisms on torsionfree  cocompact Fuchsian groups grow at most linearly in $\ell_\gamma$. 
\begin{proposition}\label{tails}
Let $M=\G\bk\h$ be a closed hyperbolic surface, and
let $f:\Gamma\to\Z$ be a conjugacy class invariant quasimorphism with defect $D$. Then 
$$
|f(\g)|\ll (1+D)\ell_\g.
$$
\end{proposition}

\begin{proof}
Fix a finite generating set $\mathcal{S}$ for $\Gamma$, and let $d_\mathcal{S}$ denote the word metric on $\Gamma$ with respect to $\mathcal{S}$. Each $\gamma\in\Gamma$ can be written as a finite reduced word in $\mathcal{S}$ so that we have the immediate bound
\begin{align*}
|f(\gamma)| \leq \sum_{\gamma_i\in\mathcal{S}}|f(\gamma_i)| +D\ell_\mathcal{S} (\gamma)\ \leq\ (C+D)\ell_\mathcal{S}(\gamma),
\end{align*}
where $D$ is the quasimorphism defect of $f$, and $C$ is an absolute constant that depends on $f$ and $\mathcal{S}$. Since $\Gamma$ acts cocompactly on $\h$, the Milnor--\v{S}varc lemma (or, as it is sometimes called, the fundamental observation of geometric group theory, see, e.g. \cite[Theorem 8.2]{FarbMargalit}) says that the metric spaces $(\Gamma,d_\mathcal{S})$ and $(\h,d_\h)$ are quasi-isometric, with the quasi-isometry given by $\Gamma\to\h$, $\gamma\mapsto \gamma z$ for any fixed $z\in\h$. In particular, there exist further constants $A\geq1$, $B\geq0$ such that for any $z$ in a fixed fundamental domain for $\G$, we have
$
\ell_\mathcal{S}(\gamma) \leq A\cdot d_\h(\gamma z,z) + B.
$
\end{proof}

We use both Theorem \ref{Thm1} and Theorem \ref{Thm:Cauchy} to prove the following equidistribution result.

\begin{theorem}[Theorem \ref{Thm:equidistribution}]
Let $\G\in\mathcal{G}$ and consider the winding number given by  Theorem \ref{Thm1c}. Let $A\subseteq\Z$ be a set with natural density $d(A)$. The density of prime geodesics with winding number in $A$ is equal to $d(A)$. Explicitly, if $\pi_A(T)$ denotes the number of prime geodesics in $\Pi(T)$ with winding number in $A$, then
$$
\lim_{T\to\infty} \frac{\pi_A(T)}{\pi(T)}\ =\ d(A).
$$
\end{theorem}

\begin{proof}
Fix $\eps>0$. There exist positive constant $K_0$, $K_1$ such that...
\begin{itemize}
\item[(i)] ...for all $K\geq K_0$,
\begin{align*}
\left|\frac{\left|\{n\in A: |n|\leq K\}\right|}{\left|\{n\in\Z:|n|\leq K\}\right|} - d(A)\right| < \eps
\end{align*}
\item[(ii)]
...and for all $K\geq K_1$,
\begin{align*}
\int_{4\pi\tfrac{K_1}{k}}^\infty \frac{du}{\pi(1+u^2)}\ <\ \eps.
\end{align*}
\end{itemize}
Choose $T$ sufficiently large that $K_1 T>K_0$. We decompose the set $A$ as
$$
A\ =\ A_1 \cup A_2 \cup A_3\ :=\ A\cap\left(\{|n|\leq K_0\}\cup\{K_0<|n|\leq K_1 T\}\cup\{K_1 T<|n|\}\right).
$$
Using (1) and applying summation by parts twice, we have
\begin{align*}
\pi_{A_2}(T)\ &=\ \sum_{\substack{|n|\leq K_1 T\\n\in A}} 1 \cdot \pi_{K_1 T}(T)-\sum_{\substack{|n|\leq K_0\\n\in A}}1\cdot \pi_{K_0}(T) - \int_{K_0}^{K_1 T} \left(\sum_{\substack{|n|\leq t\\ n\in A}}1\right) d\pi_t(T)\\
 &=\ d(A)\left(\sum_{\substack{|n|\leq K_1 T}} 1 \cdot \pi_{K_1 T}(T)-\sum_{\substack{|n|\leq K_0}}1\cdot \pi_{K_0}(T) - \int_{K_0}^{K_1 T} \left(\sum_{\substack{|n|\leq t}}1\right) d\pi_t(T)\right)\\
 &\qquad  + O\left(\eps \pi(T)\right)\ =\ d(A)\sum_{K_0<|n|\leq K_1 T} \pi_n(T) + O\left(\eps\pi(T)\right).
\end{align*}
Then
\begin{align*}
\left|\frac{\pi_A(T)}{\pi(T)}-d(A)\right|\ &\leq\ \frac{\pi_{A_1\cup A_3}(T)}{\pi(T)} + d(A)\left|\sum_{K_0<|n|\leq K_1 T}\frac{\pi_n(T)}{\pi(T)} -1\right| +O(\eps) \\
&=\ \frac{\pi_{A_1\cup A_3}(T)}{\pi(T)} +d(A)\left( \sum_{|n|\leq K_0} \frac{\pi_n(T)}{\pi(T)} +\sum_{|n|>K_1 T} \frac{\pi_n(T)}{\pi(T)}\right) +O(\eps)\\
&\leq\ (1+d(A))\left(\sum_{|n|\leq K_0}\frac{\pi_n(T)}{\pi(x)}+\sum_{|n|>K_1 T}\frac{\pi_n(T)}{\pi(T)}\right) +O(\eps).
\end{align*}
From (\ref{IBP}) and the prime geodesic theorem, the first sum is of growth order $O(\tfrac{1}{T})$ as $T\to\infty$.  Combining Theorem \ref{Thm:Cauchy} and (2), we have that
$$
\sum_{|n|> K_1 T} \frac{\pi_n(T)}{\pi(T)} =2 \int_{4\pi\tfrac{K_1}{k}}^\infty \frac{du}{\pi(1+u^2)} + O\op{\frac{1}{T}} <\eps+O\op{\frac{1}{T}}.
$$
We conclude by letting $T\to\infty$ and choosing $\eps$ to be arbitrarily small.
\end{proof}

\thanks{
The first named author would like to thank Marc Burger and Alessandra Iozzi for first mentioning winding numbers in relation to the work of Goldstein, Jay Jorgenson for discussions, and the Hausdorff Institute for Mathematics in Bonn, where most of the final draft of this manuscript was completed. The second author would like to thank Morten Risager for advice and help. C.B. is currently supported by the Swiss National Science Foundation, Grant No.~201557.}

%\bibliographystyle{alpha}
%\bibliography{biblio}

\end{document}